\documentclass[a4paper,reqno]{amsart}

\usepackage[utf8]{inputenc}
\usepackage[T1]{fontenc}
\usepackage{lmodern}
\usepackage[draft=false,kerning=true]{microtype}
\usepackage[english]{babel}

\usepackage{suffix}
\usepackage{algorithmic}
\usepackage{algorithm}
\usepackage{amssymb, bbm}

\newtheorem{theorem}{Theorem}[section]
\newtheorem{lemma}[theorem]{Lemma}
\newtheorem{proposition}[theorem]{Proposition}

\theoremstyle{definition}
\newtheorem{definition}[theorem]{Definition}

\theoremstyle{remark}
\newtheorem{remark}[theorem]{Remark}
\newtheorem{algo}[theorem]{Algorithm}

\newcommand{\nbd}{\nobreakdash}

\newcommand{\parentheses}[4][]%
{\ifthenelse{\equal{#1}{}}{\left#2}{\csname#1\endcsname#2}%
    {#4}\ifthenelse{\equal{#1}{}}{\right#3}{\csname#1\endcsname#3}}
\newcommand{\foperator}[1]{\ensuremath{%
    \mathop{{#1}\thinspace\negthinspace}
    \mathchoice{\negthinspace}{\negthinspace}{}{}}}
\newcommand{\f}[3][]{\ensuremath{\foperator{#2}\parentheses[#1]{(}{)}{#3}}}
\WithSuffix\newcommand{\f}*[2]{\ensuremath{
    \mathop{{#1}\thinspace\negthinspace}{#2}}}

\newcommand{\setmiddlesymbol}{:}
\newcommand{\set}[1]{\ensuremath{\parentheses{\{}{\}}{#1}}}
\WithSuffix\newcommand{\set}*[2]{\ensuremath{%
    \left\{{#1}\thinspace\setmiddlesymbol\thinspace{#2}\right\}}}

\makeatletter
\newcommand\undisp[1]{\bgroup\@displayfalse #1\egroup}
\makeatother

\newcommand{\conj}[1]{\ensuremath{\overline{#1}}}

\newcommand{\Z}{\ensuremath{\mathbbm{Z}}}
\newcommand{\Q}{\ensuremath{\mathbbm{Q}}}
\newcommand{\R}{\ensuremath{\mathbbm{R}}}
\newcommand{\C}{\ensuremath{\mathbbm{C}}}

\newcommand{\cD}{\ensuremath{\mathcal{D}}}
\newcommand{\cN}{\ensuremath{\mathcal{N}}}

\newcommand{\re}[2][]{\ensuremath{\f[#1]{\operatorname{Re}}{#2}}}
\newcommand{\im}[2][]{\ensuremath{\f[#1]{\operatorname{Im}}{#2}}}
\newcommand{\abs}[2][]{\ensuremath{%
    \parentheses[#1]{\lvert}{\rvert}{#2}}}

\newcommand{\ceil}[1]{\ensuremath{\parentheses{\lceil}{\rceil}{#1}}}

\newcommand{\Ztau}{\ensuremath{\Z[\tau]}}
\newcommand{\wNAF}[1][w]{$#1$\nbd-NAF}

\newcommand{\wt}[1]{\ensuremath{\widetilde{#1}}}

\newcommand{\divides}{\ensuremath{\mathbin{\mid}}}
\newcommand{\ndivides}{\ensuremath{\mathbin{\nmid}}}

\newcommand{\normsymbol}{\ensuremath{\cN}}
\newcommand{\fieldnorm}[2][]{\ensuremath{\f[#1]{\normsymbol}{#2}}}
\newcommand{\norm}[1]{\ensuremath{\parentheses{\lVert}{\rVert}{#1}}}

\usepackage{tikz}
\usetikzlibrary{arrows, snakes, decorations.pathreplacing, spy, positioning,
  calc, fit}

\usepackage{enumitem}

\usepackage{hyperref}

\let\oldpart\part
\makeatletter
\renewcommand{\part}[1]{%
  \def\@secnumfont{\bfseries}%
  \oldpart{#1}%
  \def\@secnumfont{\mdseries}}
\makeatother

\DeclareMathOperator{\Gal}{Gal}

\newcommand{\tmod}{\!\!\mod}

\newcommand{\fp}{\mathfrak{p}}
\newcommand{\fpbar}{\bar{\mathfrak{p}}}
\newcommand{\cJ}{\mathcal{J}}
\newcommand{\cJbar}{\bar{\mathcal{J}}}

\numberwithin{equation}{section}
\numberwithin{figure}{section}
\numberwithin{table}{section}
\numberwithin{algorithm}{section}

\begin{document}

\newcommand{\qupper}{500}


\title[Non-Minimality of the Width-$w$ Non-adjacent Form]{%
  Non-Minimality\\ of the Width-$w$ Non-adjacent Form\\
  in Conjunction with\\
  Trace One $\tau$-adic Digit Expansions and\\
  Koblitz Curves in Characteristic Two}

\author{Daniel Krenn}

\thanks{Daniel Krenn is supported 
  by the Austrian Science Fund (FWF): I1136,
  by the Austrian Science Fund (FWF): P24644-N26, and
  by the Austrian Science Fund (FWF): W1230, 
  Doctoral Program ``Discrete Mathematics''.}

\address{\parbox{12cm}{%
    Daniel Krenn \\
    Institut f\"ur Mathematik \\
    Alpen-Adria-Universit\"at Klagenfurt \\
    Universit\"atsstra\ss e 65--67, 9020 Klagenfurt am W\"orthersee, Austria\\}}

\email{\href{mailto:math@danielkrenn.at}{math@danielkrenn.at} \textit{or}
  \href{mailto:daniel.krenn@aau.at}{daniel.krenn@aau.at}}

\author{Volker Ziegler}

\thanks{Volker Ziegler is supported
  by the Austrian Science Fund (FWF): P24801.}

\address{\parbox{12cm}{%
    Volker Ziegler \\
    Fachbereich für Mathematik \\
    University of Salzburg  \\
    Hellbrunnerstrasse 34, A-5020 Salzburg, Austria \\}} 

\email{\href{mailto:volker.ziegler@sbg.ac.at}{volker.ziegler@sbg.ac.at}}

\keywords{$\tau$-adic expansions, redundant digit
  sets, elliptic curve cryptography, Koblitz curves, Frobenius endomorphism,
  scalar multiplication, Hamming weight, linear forms in logarithms, geometry
  of numbers, Baker--Davenport method, continued fractions}

\subjclass[2010]{%
11A63, 
11Y50, 
11D75. 
}


\begin{abstract}
  This article deals with redundant digit expansions with an imaginary
  quadratic algebraic integer with trace $\pm 1$ as base and a minimal norm
  representatives digit set. For $w\geq 2$ it is shown that the width-$w$
  non-adjacent form is not an optimal expansion, meaning that it does not
  minimize the (Hamming-)weight among all possible expansions with the same
  digit set. One main part of the proof uses tools from Diophantine analysis,
  namely the theory of linear forms in logarithms and the
  Baker--Davenport reduction method.
\end{abstract}


\maketitle


\part{The Beginning}
\label{part:I}


\section{Introduction}
\label{sec:intro}


Let $\tau$ be an (imaginary quadratic) algebraic integer and $\cD$ a finite
subset of $\Ztau$ including zero. Choosing the digit set $\cD$ properly, we can represent
$z\in\Ztau$ by a finite sum
\begin{equation*}
  \sum_{\ell=0}^{L-1} \sigma_\ell \tau^\ell,
\end{equation*}
where the digits $\sigma_\ell$ lie in $\cD$. We call this representation a digit
expansion of $z$. Using a redundant digit set $\cD$, i.e., taking more digits
than needed to represent all elements of $\Ztau$, each element can be written
in different ways. Of particular interest are expansions
which have the lowest number of nonzero digits. We call those expansions
\emph{optimal} or \emph{minimal expansions}.

The motivation looking at such expansions comes from elliptic curve
cryptography. There the scalar multiplication of a point on the curve is a
crucial operation and has to be performed as efficiently as possible. The
standard double-and-add algorithm can be extended by windowing methods, see for
example \cite{Blake-Seroussi-Smart:1999, Gordon:1998,
  Miyaji-Ono-Cohen:1997:effic, Solinas:1997:improved-algorithm,
  Solinas:2000:effic-koblit}. Translating this into the language of digit
expansions means the usage of redundant digit expansions with base~$2$. However,
using special elliptic curves, for example Koblitz curves, see
\cite{Koblitz:1992:cm, Koblitz:1998:ellip-curve,
  Solinas:1997:improved-algorithm, Solinas:2000:effic-koblit}, the
``expensive'' doublings can be replaced by the ``cheap'' application of the
Frobenius endomorphism over finite fields. In the world of digit expansions
this means taking an imaginary quadratic algebraic integer as base. This leaves 
us with the additions of points of the elliptic curve as an ``expensive''
operation. The number of such additions is basically the number of
nonzero digits in an expansion. Therefore minimizing
this number is an important goal.

We are now going back to expansions with a low number of nonzero digits.
Let the parameter $w\geq2$ be an integer.
Then one special expansion is the width\nbd-$w$ non-adjacent form,
where in each block of width $w$ at most one digit is not equal to zero, see
Reitwiesner~\cite{Reitwiesner:1960} who introduced this notion for $w=2$ and
others including Muir and Stinson~\cite{muirstinson:minimality} and
Solinas~\cite{Solinas:1997:improved-algorithm, Solinas:2000:effic-koblit}. It
will be abbreviated by \wNAF{}.  This expansion contains, by construction, only
few nonzero digits. When we use a digit set consisting of zero and of
representatives with minimal norm of the residue classes modulo $\tau^w$
excluding those which are divisible by $\tau$, then the \wNAF{}-expansion is optimal
(minimal) in a lot of cases. For example, using an integer (absolute value at
least $2$) as base~$\tau$, the \wNAF{} is a minimal/optimal expansion, see
Reitwiesner~\cite{Reitwiesner:1960}, Jedwab and
Mitchell~\cite{Jedwab-Mitchell:1989}, Gordon~\cite{Gordon:1998},
Avanzi~\cite{avanzi:mywnaf}, Muir and Stinson~\cite{muirstinson:minimality},
and Phillips and Burgess~\cite{Phillips-Burgess:2004:minim-weigh}. As a digit
set it contains in these cases zero and all integers with absolute value smaller than
$\frac12\abs{\tau}^w$ and not divisible by $\tau$.

A general criterion for optimality of the \wNAF{}-expanions can be found in
Heuberger and Krenn~\cite{Heuberger-Krenn:2013:wnafs-optimality}: The \wNAF{}
of each element is optimal, if expansions of weight~$2$ are optimal. This is
especially useful if the digit set has some underlying geometric properties
as it is the case for a minimal norm representatives digit set. In Heuberger
and Krenn~\cite{Heuberger-Krenn:2013:general-wnaf} an optimality result for a
general algebraic integer base is given. A refinement
of this general criterion in the imaginary quadratic case is stated in
\cite{Heuberger-Krenn:2013:wnafs-optimality}.
For $\tau$ being imaginary quadratic and a zero
of $\tau^2-p\tau+q$, the main result is that optimality follows if
$\abs{p}\geq3$ and $w\geq4$. Further, there are conditions given for $w=2$ and
$w=3$. In the cases $p=\pm2$ and $q=2$ the \wNAF{}-expansion is optimal for odd
$w$ and non-optimal for even $w$. Moreover, non-optimality was also shown when
$p=0$ and $w$ is odd.

In this article we are interested in the case when $p\in\set{-1,1}$. Note that
the case $q=2$ is related to Koblitz curves in characteristic $2$,
see Koblitz~\cite{Koblitz:1992:cm}, Meier and
Staffelbach~\cite{Meier-Staffelbach:1993:effic}, and
Solinas~\cite{Solinas:1997:improved-algorithm, Solinas:2000:effic-koblit}. A few
results are already known: If $w=2$ or $w=3$ optimality was shown in
Avanzi, Heuberger and
Prodinger~\cite{Avanzi-Heuberger-Prodinger:2006:minim-hammin,
  Avanzi-Heuberger-Prodinger:2006:scalar-multip-koblit-curves} (see also
Gordon~\cite{Gordon:1998} for $w=2$). In contrast, for $w\in\set{4,5,6}$, the
\wNAF{}-expansion is not optimal anymore. This was shown in
Heuberger~\cite{Heuberger:2010:nonoptimality}. These results rely on transducer
automata rewriting arbitrary expansions (with given base and digit set) to
a \wNAF{}-expansion and on a search of cycles of negative ``weight''.

Experimental results checking the above criterion by symbolic
calculations indicate that the \wNAF{} is non-optimal for $w\geq4$ and,
moreover, non-optimal for all $w\geq2$ when $q\geq3$, see Heuberger and
Krenn~\cite{Heuberger-Krenn:2013:wnafs-optimality}. The main result presented
in this work---see the next section for a precise statement---proves
this conjecture for $q\leq \qupper$.


\section{Expansions, the Results and an Overview}
\label{sec:overview}


We use this section to present our main theorem and to give an overview of the
different methods used during its proof. We start by explaining what we mean by
optimal (or minimal) expansions.


\begin{definition}\label{def:expansion}
  Let $\tau$ be an algebraic integer, and suppose we have a set $\cD$ (called
  \emph{digit set}) with $\cD\subseteq\Ztau$ such that $0\in\cD$. Let $w$
  be an integer with $w\geq2$ (called \emph{window size}).
  \begin{enumerate}
  \item The finite sum
    \begin{equation*}
      z = \sum_{\ell=0}^{L-1} \sigma_\ell \tau^\ell
    \end{equation*}
    with a positive integer~$L$ and $\sigma_\ell\in\cD$ for all $\ell$ is called
    a \emph{digit expansion of~$z$} with base~$\tau$.

  \item We call the expansion defined above a \emph{width\nbd-$w$ non-adjacent
      form} (abbreviated by \emph{\wNAF{}}) if for
    $\ell\in\set{0,1,\dots,L-w}$ each of the words
    \begin{equation*}
      \sigma_{\ell}\sigma_{\ell+1}\dots\sigma_{\ell+w-1}
    \end{equation*}
    contains at most one nonzero digit.

  \item The number of nonzero digits is called the \emph{(Hamming-)weight} of
    the expansion.

  \item Suppose we have an expansion of~$z$ with weight $W$. We call this
    expansion \emph{optimal} or \emph{minimal} if each expansion of~$z$ (with
    digits out of $\cD$) has a weight which is at least $W$.

  \item The \wNAF{} expansion is said to be \emph{optimal} or
    \emph{minimal} (with respect to $\tau$ and to $\cD$) if the \wNAF{}
    of each element of $\Ztau$ is minimal.

  \end{enumerate}
\end{definition}


We will skip ``with respect to $\tau$ and to $\cD$'' in the previous
definitions if this is clear from the context (and in our cases it will always
be the base $\tau$ and the minimal norm digit set~$\cD$).


Before we are able to state our main theorems, we have to specify the digit
set~$\cD$. For a parameter~$w$ (the window size) we assume that $0$ is a digit
and that we take a representative of minimal norm out of each residue class
modulo $\tau^w$ which is not divisible by the base~$\tau$. We call such a digit
set a \emph{minimal norm representative digit set modulo $\tau^w$}, see
Section~\ref{sec:digit-sets}, in particular
Definition~\ref{def:min-norm-digit-set}, for details. 


\begin{remark}
  If $\tau$ is an imaginary quadratic algebraic integer (as we use it here in
  this article) and $\cD$ a minimal norm representative digit set modulo
  $\tau^w$, then each element of $\Ztau$ admits a unique \wNAF{} expansion,
  see Heuberger and Krenn~\cite{Heuberger-Krenn:2013:wnaf-analysis}.
\end{remark}


Now it is time to state our main results.


\begin{theorem}\label{thm:non-opt-large-w}
  Let $q$ be an integer with $q\geq2$ and let $p\in\set{-1,1}$. Let $\tau$ be a
  root of $X^2-pX+q$ and $\cD$ a minimal norm representative digit set modulo
  $\tau^w$. Then there exists an effectively computable bound~$w_q$ such that
  for all $w\geq w_q$ the width\nbd-$w$ non-adjacent form expansion is not
  optimal with respect to $\tau$ and to $\cD$. In particular, we may
  choose\footnote{The explicit bounds $w_q$
    (Theorem~\ref{thm:non-opt-large-w}, ``in particular''-part) are rough
    estimates. For a particular $q$, better bounds can be
    computed, which is done throughout this article.}
  \begin{align*}
    w_q&=8.68\cdot 10^{15} \log q \log\log q &&\text{if $q\geq 13$}
    \intertext{and}
    w_q&=1.973 \cdot 10^{16} &&\text{if $q\in\set{2, 3, \dots, 12}$}.
  \end{align*}
\end{theorem}


It turns out that the bounds are rather huge (Section~\ref{sec:bounds-w}).
However, for small~$q$ we can reduce this bound
dramatically (for example from $w_2=8.596\cdot 10^{15}$ to $\wt w_2=140$)
and get the following much stronger result.


\begin{theorem}\label{thm:main-non-opt}
  Let $q$ and $w$ be integers with
  \begin{itemize}
  \item either\footnote{The upper bound~$q\leq\qupper$ in
      Theorem~\ref{thm:main-non-opt}
      is determined by an extensive computation. Details can be found at
      \url{http://www.danielkrenn.at/koblitz2-non-optimal}, in particular file
      \href{http://www.danielkrenn.at/downloads/koblitz2-non-optimal/result_overview}{\texttt{result\_overview}}.}
    $2\leq q \leq \qupper$ and $w\geq2$
  \item or $q\geq2$ and $w\in\set{2,3}$,
  \end{itemize}
  and let $p\in\set{-1,1}$. Let $\tau$ be a root of $X^2-pX+q$ and $\cD$ a
  minimal norm representative digit set modulo $\tau^w$. Then the width\nbd-$w$
  non-adjacent form expansion with respect to $\tau$ and to $\cD$ is
  optimal if and only if $q=2$ and $w\in\set{2,3}$.
\end{theorem}


Formulated differently, this means that the \wNAF{} is \textbf{not}
minimal/optimal for all the given parameter configurations except for the four
cases with $w\in\set{2,3}$, $p\in\set{-1,1}$ and $q=2$.


The main part of the proof of Theorem~\ref{thm:main-non-opt} deals with an
algorithm which takes $q$ (and $p$) as input and outputs a list of values for $w$
for which no counterexample to the minimality of the \wNAF{} was found. Let us
also formulate this as a proposition.


\begin{proposition}\label{pro:main-algo}
  Let $q$ be an integer with $q\geq2$ and $p\in\set{-1,1}$. Let $\tau$ be a
  root of $X^2-pX+q$ and $\cD$ a minimal norm representative digit set modulo
  $\tau^w$. Then there is an algorithm which tests non-optimality of the
  width\nbd-$w$ non-adjacent form expansion for all $w\geq2$.
\end{proposition}


This algorithm grows out of an intuition on how a counterexample to minimality
of the \wNAF{} are constructed. To do so, we have to find certain lattice
point configurations located near the boundary of the digit set. This is
described in general in the overview (Section~\ref{sec:digits-overview}) of
Part~\ref{sec:digits} and with more details and very specific for our situation
in Section~\ref{sec:strategy}. This leaves us to find a lattice point located
in some rectangle which additionally avoids some smaller lattice. 

The whole Part~\ref{sec:dioph} deals with this problem of finding a suitable
lattice point inside the given rectangle, which is precisely formulated as
Proposition~\ref{pro:find-lattice-point}. Using the theory of the geometry of numbers
allows us to construct such a lattice point, but unfortunately not ``for
free''; we have to ensure that there is no lattice point in some smaller
rectangle. This problem can be reformulated as an inequality (namely
Inequality~\eqref{eq:Main_Log_Form}) and we have to show that it does not have
any integer solutions.

Dealing with the solutions of Inequality~\eqref{eq:Main_Log_Form} is a task of
Diophantine analysis. In particular and because of the structure
of~\eqref{eq:Main_Log_Form-a} we use the theory of linear forms in
logarithms. This provides, for given $q$, a rather huge bound on $w$ (due to
Matveev~\cite{Matveev:2000:explicit-lower-bound}), see
Section~\ref{sec:bounds-w} for details. From this
Theorem~\ref{thm:non-opt-large-w} can be proven. However, using the convergents
of continued fractions we are able to reduce this bound
significantly. Therefore, we are able to check all the remaining $w$
directly. In particular, we use a variant of the Baker--Davenport
method~\cite{Baker-Davenport:1969}, which is described in
Section~\ref{sec:smaller-bound}.

Let us close this overview of the second part with the following: In
Section~\ref{sec:testing-fixed-w} we give some remarks on how to test
Inequality~\eqref{eq:Main_Log_Form} directly. The actual algorithm is stated in
Section~\ref{sec:algorithm}.

In Part~\ref{sec:digits} digits come into play and the counterexamples to
minimality of the \wNAF{} are constructed. Section~\ref{sec:non-opt-w} explains
this directly for some values of~$w$ (but arbitrary~$q$), whereas the remaining
sections deal with the construction using the result of
Part~\ref{sec:dioph}. In particular, this gives us a minimal non-\wNAF{}
expansion, whose most significant digit is perturbated a little bit
(Section~\ref{sec:small-change-Delta}).
This is compensated by a large change in
the least significant digit, see Section~\ref{sec:lsd}. In the
Sections~\ref{sec:strategy} and~\ref{sec:proof-main} all results are glued
together and the actual counterexamples are constructed (thereby proving
Theorem~\ref{thm:main-non-opt}). The actual algorithm is
implemented\footnote{See
  \url{http://www.danielkrenn.at/koblitz2-non-optimal} for the code.}
in SageMath~\cite{SageMath:2016:7.0}.

We are now finished with the introductory overview and will start with two
preparatory sections.


\section{The Set-Up}
\label{sec:setup}


This section is to state a couple of definitions used throughout this work and
to fix some notations.


\begin{itemize}


\item Let $q$ be an integer with $q\geq2$. We call $q$ the \emph{norm of our
    base} (for what we mean by ``base'', have a look below).


\item Let $p\in\set{-1,1}$. We call $p$ the \emph{trace of our base}.


\item Let $\tau$ be a zero of $X^2 - p X + q$, more precisely, we take
  \begin{equation*}
    \tau = \frac{p}{2} + i\sqrt{q-\frac14}.
  \end{equation*}
  We call $\tau$ the \emph{base of our expansions}.

  Note that the case $\tau = \tfrac12p - i\sqrt{q-\tfrac14}$ (i.e., taking the
  negative square root) is, by conjugation, ``equivalent'' to our set-up. This
  shall mean by constructing a counterexample to minimality of the \wNAF{} in
  one case (sign of the square root) and by conjugating everything, we obtain a
  counterexample for the other sign of the square root.


\item Our expansions live in the lattice
  \begin{equation*}
    \Ztau = \set*{a+b\tau}{\text{$a\in\Z$, $b\in\Z$}}.
  \end{equation*}
  See also Section~\ref{sec:lattices} for the other lattices used.


\item We set
  \begin{equation*}
    V = \set*{z\in\C}{ \forall y\in\Ztau \colon \abs{z} \leq \abs{z-y}}
  \end{equation*}
  and call it the \emph{Voronoi cell of $0$} corresponding to the set
  $\Ztau$. An example of this Voronoi cell in a lattice $\Ztau$ is shown in
  Figure~\ref{fig:voronoi}.

  \begin{figure}
    \centering
    \includegraphics{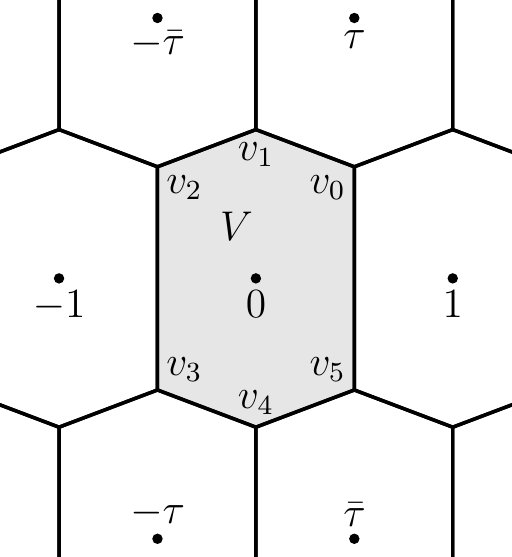}
    \caption[Voronoi cell $V$ for $0$]{Voronoi cell $V$ of $0$ corresponding
      to the set $\Ztau$ with $\tau = \frac{1}{2}+\frac{i}{2}\sqrt{7}$ (i.e.,
      $p=1$, $q=2$).}
    \label{fig:voronoi}
  \end{figure}


\item The \emph{vertices of $V$} are
  \begin{align*}
    v_0 &= \frac{p}{2} + \frac{i}{2\im{\tau}} \left(\im{\tau}^2-\frac14\right)
    = \frac{p}{2} + \frac{i}{\sqrt{4q-1}} \left(q-\frac12\right)\!, \\
    v_1 &= \frac{i}{2\im{\tau}} \left(\im{\tau}^2+\frac14\right)
    = \frac{i}{\sqrt{4q-1}} q,  \\
    v_2 &= -\frac{p}{2} + \frac{i}{2\im{\tau}} \left(\im{\tau}^2-\frac14\right)
    = -\frac{p}{2} + \frac{i}{\sqrt{4q-1}} \left(q-\frac12\right)\!,
  \end{align*}
  and $v_3=-v_0$, $v_4=-v_1$ and $v_5=-v_2$, see Heuberger and
  Krenn~\cite{Heuberger-Krenn:2013:wnaf-analysis}.


\item Let $w$ be an integer with $w\geq2$. We call $w$ the \emph{window size}
  of our expansions, see also Definition~\ref{def:expansion}.


\item Let
  \begin{equation*}
    d = \frac{v_1-v_0}{\abs{v_1-v_0}}
  \end{equation*}
  be the direction from $v_0$ to $v_1$, and note that we have
  \begin{equation*}
    d = pi\frac{\tau}{\sqrt{q}}.
  \end{equation*}
  Set $d_w = d (\tau/\sqrt{q})^w$. See also Figure~\ref{fig:rectangle-Rw}.


\item Let
  \begin{equation*}
    s = \sqrt{\frac{q-1/4}{q+2}}
  \end{equation*}
  be the \emph{height of the rectangle} defined below. See also
  Figure~\ref{fig:rectangle-Rw}.


\item Define the open \emph{rectangle~$R_w$} with vertices
  \begin{itemize}
  \item $\tau^w v_0 + d_w\sqrt{q}$,
  \item $\tau^w v_1 - d_w\sqrt{q}$,
  \item $\tau^w v_1 - d_w\sqrt{q} - pid_ws$, and
  \item $\tau^w v_0 + d_w\sqrt{q} - pid_ws$.
  \end{itemize}
  An example of this rectangle is shown in Figure~\ref{fig:rectangle-Rw}.

  \begin{figure}
    \centering
    \includegraphics{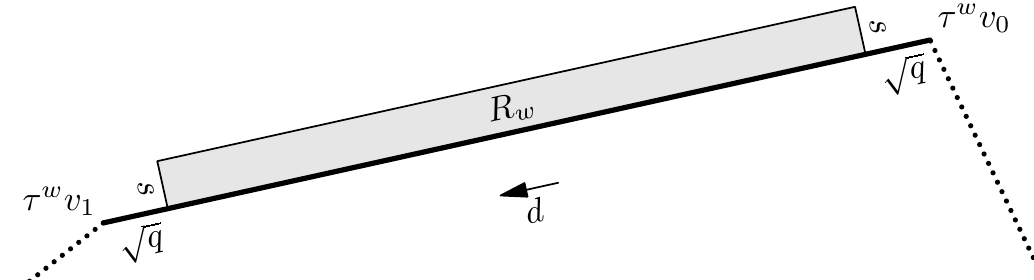}
    \caption{Rectangle~$R_5$ for $q=5$ and $p=1$.}
    \label{fig:rectangle-Rw}
  \end{figure}
  
  Note that one side length of the rectangle $R_w$ is
  \begin{equation}\label{eq:width-Rw}
    \abs{\tau^w v_1 - \tau^w v_0} - 2\sqrt{q} 
    = \sqrt{q}^w \frac12 \sqrt{1+\frac{1}{4q-1}} - 2\sqrt{q}
    = \sqrt{\frac{q^{w+1}}{4q-1}} - 2\sqrt{q}
  \end{equation}
  and the other is~$s$.


\end{itemize}


We finish this section with a couple of remarks.

\begin{remark}\label{rem:Rw-divide}
  The location of the rectangle~$R_w$ (in relation to the scaled Voronoi cell
  $\tau^w V$) is in such a way that $\tau^{-1} R_w$ has empty intersection with
  $\tau^w V$. When constructing the actual counterexample in
  Part~\ref{sec:digits}, this will guarantee us that an element of $R_w$ does
  not become a digit (during division by $\tau$).
\end{remark}

\begin{remark}\label{rem:Rw-width}
  Note that the rectangle~$R_w$ is well defined (has positive area) if and only
  if $w>\log(16q-4) / \log q$, which follows from positivity
  of~\eqref{eq:width-Rw}. Therefore, for $q=2$ we have $w\geq5$, for $q=3$ we
  have $w\geq4$, for $4\leq q\leq 15$ we have $w\geq3$ and for $q\geq 16$ we
  have $w\geq2$.
\end{remark}


\section{Lattices}
\label{sec:lattices}


As mentioned above, our digit expansions live in the lattice
\begin{equation*}
  \Ztau = \set*{a+b\tau}{\text{$a\in\Z$, $b\in\Z$}} = \langle 1, \tau \rangle.
\end{equation*}
It will become handy to define a few other (related) lattices. Our first one is
\begin{equation*}
  \Lambda_\tau=\langle \tau, q-p\tau\rangle=\langle \tau, \tau^2 \rangle,
\end{equation*}
where we interpret the complex plane embedded into $\R^2$ in the usual way. This
lattice is used during the construction of our counterexamples, since we need
points divisible by $\tau$ there. Further, we also work with the smaller
lattice
\begin{equation*}
  \Lambda_{\tau^2}
  =\langle \tau^2, \tau^3 \rangle
  = \langle q\tau, \tau^2 \rangle
  \subseteq \Lambda_\tau,
\end{equation*}
since in view of Proposition~\ref{pro:find-lattice-point} we want to avoid this
lattice.

Moreover, let us note that the middle point of the lower long side of the
rectangle~$R_w$ is
\begin{equation*}
\frac{v_0+v_1}2 \tau^w=\frac{\tau^{w+1}}2.
\end{equation*}
In general this is not a point of the lattice~$\Lambda_\tau$ but of the lattice
\begin{equation*}
  \Lambda_{\tau/2}
  = \frac 12 \Lambda_\tau
  = \langle \frac{\tau}{2}, \frac{\tau^2}{2} \rangle
  \supseteq \Lambda_\tau.
\end{equation*}
This is the reason why we will work mainly in the larger
lattice~$\Lambda_{\tau/2}$.


We need some basic properties of the lattices above. Let us start with
$\Lambda_\tau$ and some divisibility conditions for its elements.


\begin{lemma}\label{lem:structure-lattice}
  The elements of $\Ztau$ divisible by $\tau$ once (i.e., divisible by~$\tau$
  but not by~$\tau^2$) are exactly the elements
  \begin{equation*}
    a\tau + b(q-p\tau)
  \end{equation*}
  with $a\in\Z$ but $q\ndivides a$ and with $b\in\Z$.

  Moreover, if $z\in\Ztau$ is a lattice point divisible by $\tau$, then $z-1$, $z+1$,
  $z+\tau-p$ and $z-\tau+p$ are not divisible by $\tau$.
\end{lemma}

\begin{figure}
  \centering
  \includegraphics[width=\linewidth]{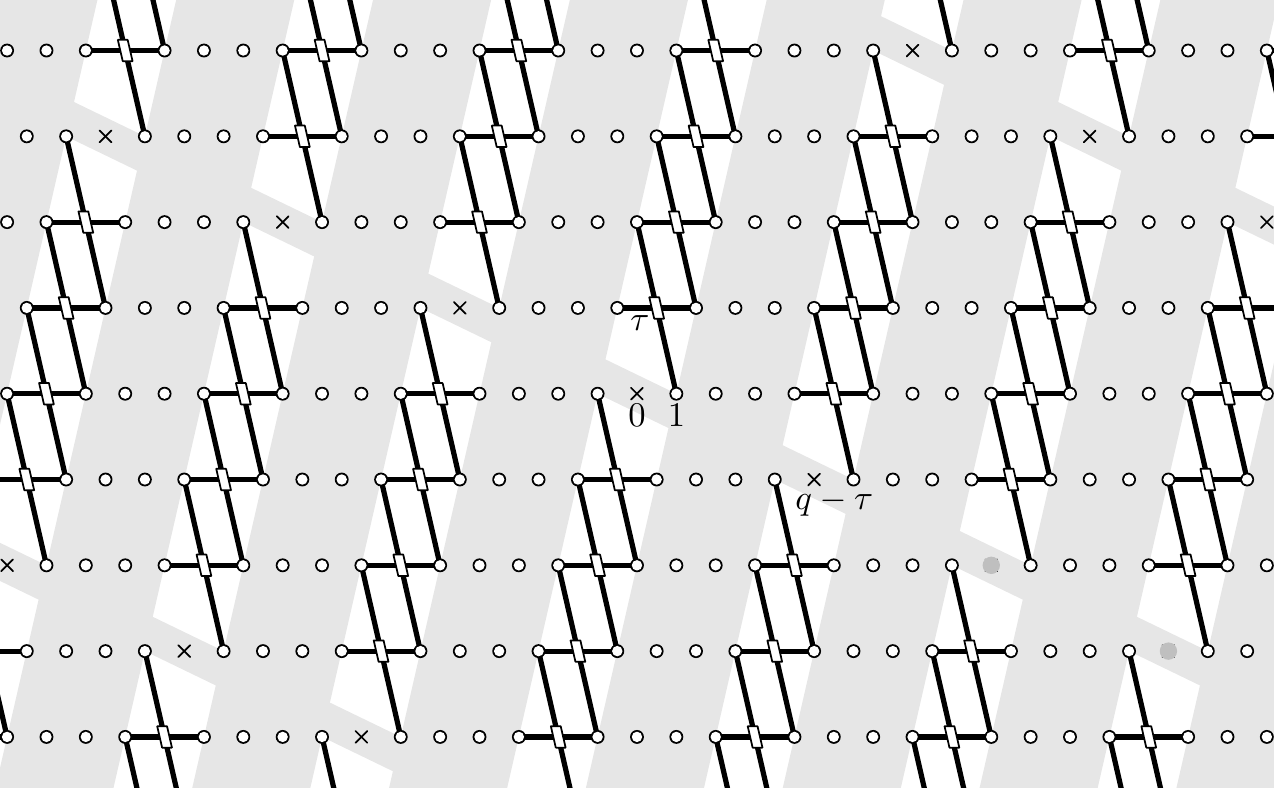}
  \caption{Lattice $\Ztau$. Points marked with a circle are not divisible
    by~$\tau$, points marked with a rectangle are divisible exactly once
    by~$\tau$, points marked with a cross are divisible by~$\tau^2$. There are
    lines from each rectangle (divisible exactly once by~$\tau$) to its
    neighboring circles (not divisible by~$\tau$)}.
  \label{fig:lattice}
\end{figure}

The structures described above can be found in Figure~\ref{fig:lattice}.

\begin{proof}[Proof of Lemma~\ref{lem:structure-lattice}]
  An element $c+d\tau \in \Ztau$ is divisible by~$\tau$ if and only if
  $q\divides c$. Therefore, if $z\in\Ztau$ is divisible by $\tau$, then
  $z-1$, $z+1$, $z+\tau-p$ and $z-\tau+p$ are not divisible by $\tau$, which
  proves the second part of the lemma.
  As $\tau\ndivides a - b\tau$ if and only if
  $q\ndivides a$, multiplying by $\tau$ yields the first part of the lemma.
\end{proof}


For analyzing the lattice $\Lambda_{\tau/2}$ and the scaled Voronoi cell $\tau^w
V$, it is important to know some arithmetic properties of $\tau$. We show the
following lemma to get some insights.


\begin{lemma}\label{lem:properties_tau}
  The algebraic integer $\tau$ satisfies the following properties.
  \begin{enumerate}
  \item Every prime $\ell\in\Z$ with $\ell\divides q$ splits in $\Z[\tau]$.
  \item If $a$ and $b$ are integers with $a+b\tau=\tau^w$ or
    $a+b\tau=\bar\tau^w$, then $\gcd(a,b)=1$.
  \end{enumerate}
\end{lemma}


Note that $\Ztau$ is the maximal order of~$\Q(\tau)$.


\begin{proof}[Proof of Lemma~\ref{lem:properties_tau}]
  The first statement is a direct consequence from algebraic number theory, in
  particular the splitting of a prime $\ell$ in $\Z[\tau]$ is described by the
  factorization of the minimal polynomial of $\tau$ mod $\ell$ (for example, see
  Theorem~2 in Chapter~11 of Ribenboim~\cite{Ribenboim:2001}). Indeed we have
  \begin{equation*}
    X^2-pX+q \equiv X^2\pm X\equiv X(X\pm 1) \mod \ell,
  \end{equation*}
  hence all primes $\ell\divides q$ split completely in $\Z[\tau]$.

  The second statement is trivial for $w\in\set{0,1}$ (note that we use
  $w\geq2$ throughout this paper anyway). Suppose $a$ and $b$ have a common
  prime factor $\ell$, then $\ell$ also has to divide
  $\fieldnorm{\bar\tau}^w=\fieldnorm{\tau}^w=q^w$ (where $\fieldnorm{\tau}$
  denotes the norm of $\tau$). Thus $\ell\divides q$. By using part (1) of this
  lemma we have $(\ell)=\fp \fpbar$ as ideals over
  $\Ztau$.

  Let us assume for a moment that both $\fp$ and $\fpbar$ divide $(\tau)$, then
  also both $\fp$ and $\fpbar$ divide $(\bar\tau)$, i.e.,
  $\tau,\bar\tau \in \fp \fpbar=(\ell)$. But this yields
  $\ell\divides \tau+\bar\tau=p$, a contradiction. Therefore let us assume now
  that $\fp \divides (\tau)$ and $\fpbar\ndivides (\tau)$. Since
  by assumption $a$ and $b$ have the common factor $\ell$, they also have the
  common factor $\fpbar$ form the ideal point of view, hence
  $\fpbar \divides (\tau)^w$ if $a+b\tau=\tau^w$, a contradiction to
  the previous discussion. Similarly, we get the contradiction $\fp \divides (\bar\tau)^w$ for the case $a+b\tau=\bar\tau^w$.
\end{proof}


It is also important to know that no lattice points are on the ``lower'' edge
of $R_w$. This result is also used to show the uniqueness of the digit set, see
Proposition~\ref{pro:digit-set-unique}.


\begin{lemma}\label{lem:lower_line}
  The following two statements hold.
  \begin{enumerate}
  \item The only lattice point in $\Lambda_{\tau/2}$ lying on the line segment
    joining the points $v_0\tau^w$ and $v_1\tau^w$ is $\frac12\tau^{w+1}$.
  \item The boundary of $\tau^w V$ has empty intersection with the
    lattice~$\Ztau$.
  \end{enumerate}
\end{lemma}


\begin{proof}
  We start by showing that there are no lattice points of $\Lambda_{1/2}=
  \langle\frac12,\frac\tau2\rangle$ except $\frac12\tau^w$ on the line going
  through $v_5\tau^w$ and $v_0\tau^w$. Every point on this line can be written
  as $\frac12 \tau^w + t i\tau^w$ with some real parameter~$t$. Let us assume
  we have a point $\lambda=\frac12\tau^w + a\frac12 + b\frac\tau2$ with $a$,
  $b\in\Z$ on this line. Furthermore, we may assume that $\gcd(a,b)=1$ (as a
  minimality condition). We deduce $2ti\tau^w=a+b\tau$.

  If $i\in\Q(\tau)$, then $4q-1$ is a perfect square which is absurd since
  $4q-1\equiv -1 \mod 4$. Let us write $\tau^w=a_w+b_w\tau$, then
  Lemma~\ref{lem:properties_tau} yields $\gcd(a_w,b_w)=1$. Thus, and since
  $i\not\in\Q(\tau)$, the only possible values for $t$ are $\pm \frac 12$.
  Indeed $2ti\tau^w$ is an algebraic integer and therefore we have $2t\in\Z[\tau]\cap \R=\Z$
  with $t\neq 0$ and $|t|<1$.
  Now, we obtain a contradiction, since $\lambda=\frac12(1\pm i)\tau^w$ is not in
  $\Lambda_{1/2}$.

  The results are now obtained as follows. Multiplying everything by $\tau$
  yields the first statement of the lemma. Starting with $\frac12 \tau^{w+1} +
  t i\tau^{w+1}$ yields that there are no points from $\Lambda_{1/2}$ except
  $\frac12\tau^{w+1}$ on the line joining $v_0\tau^w$ and $v_1\tau^w$. Note
  that this differs from the first statement of the lemma by the different
  lattice. The result for the third line (from $v_1\tau^w$ to $v_2\tau^w$)
  follows by taking conjugation and Lemma~\ref{lem:properties_tau} with
  $\bar\tau^w = a_w+b_w\tau$. The remaining three sides of $\tau^w V$ follow by
  mirroring.
\end{proof}


We are also interested in the shortest vector in the lattice $\Lambda_{\tau/2}$.


\begin{lemma}\label{lem:shortest_vector}
  The shortest nonzero vectors in the lattice $\Lambda_{\tau/2}$ are
  $\pm\tau/2$.
\end{lemma}


\begin{proof}
  First, let us note that
  $\abs[big]{\frac{\tau}{2}}=\frac{\sqrt{q}}{2}$. To find the shortest
  nonzero vector $a\frac{\tau}{2}+b\frac{\tau^2}{2} \in
  \Lambda_{\tau/2}$ we have to find all integers $a$ and $b$ such that
  $\abs[big]{a\frac{\tau}2+b\frac{\tau^2}2} \leq \frac{\sqrt{q}}{2}$. In
  particular we have to solve the inequality
  \begin{equation*}
    \abs{a+b\tau}^2 =
    \Bigl(a+b\frac{p}{2}\Bigr)^2 + b^2\Bigl(q-\frac 14\Bigr)
    \leq 1.
  \end{equation*}
  Obviously, if $b\neq 0$, then this inequality cannot be
  satisfied. Thus we may assume that $b=0$. We obtain $a^2\leq
  1$, and the result follows.
\end{proof}


\part{The Diophantine Part}
\label{sec:dioph}


\section{Overview}
\label{sec:dioph-overview}


In this part of the article we show that the following proposition holds.


\begin{proposition}\label{pro:find-lattice-point}
  We use the set-up described in Section~\ref{sec:setup} with the following
  restrictions. Suppose we are in one of the cases
  \begin{itemize}
  \item $q\in\set{2,4}$ and $w\geq7$,
  \item $q\in\set{3}$ and $w\geq5$ or
  \item $5\leq q\leq \qupper$ and $w\geq4$.
  \end{itemize}
  Then there exists a lattice point
  \begin{equation*}
    a \tau + b (q-p\tau)
  \end{equation*}
  with $a\in\Z$, $q\ndivides a$ and with $b\in\Z$ in the (open) rectangle~$R_w$.
\end{proposition}


Since the proof of this proposition is long and technical, we start with an
overview. In a nutshell, for fixed~$q$, we can reduce the problem to checking
only finitely many configurations~$w$. Therefore the testing is possible
algorithmically.

Let us look at an outline of the ideas used during the proof a bit. The
resulting algorithm takes as input parameters $q$ and $p$ and returns a list of
values for $w$ which have to be investigated by other methods (i.e., yet no
lattice point was found for these~$w$). The details can be found in the last
section of this part, Section~\ref{sec:algorithm}. In order to make this
algorithm work, we have to check Proposition~\ref{pro:find-lattice-point} for all but
finitely many cases.

One major step to tackle this lemma is to reduce the lattice point problem into
a Diophantine approximation problem. We show this in
Section~\ref{sec:Geometry_of_Numbers}. The existence result of the lattice
points there is based on the theory of the geometry of numbers. More precisely,
this gives us two lattice points (see Lemma~\ref{lem:exist-lambda}) out of which we can
construct a lattice point avoiding a smaller lattice (as it is required by
Proposition~\ref{pro:find-lattice-point}), see
Lemma~\ref{lem:conditional_existence}. But, to make this work, we have to use
linear independence of the two points, which is the challenging part during the
proof.

We can reformulate this linear independence problem geometrically, which leaves
us to show that there are no lattice points inside a certain smaller
rectangle. To solve it, we bring this Diophantine approximation problem into a
favorable form, which leaves us to show that the inequalities
\begin{subequations}\label{eq:Main_Log_Form}
  \begin{equation}\label{eq:Main_Log_Form-a}
    \abs{\log\left(\frac{a+b\tau}{\abs{a+b\tau}}\right)
      - w\log\left(\frac{\tau}{\abs{\tau}}\right)
      + k\frac{i\pi}2}
    < \chi q^{2-w/2}
  \end{equation}
  with $\chi=9$ and
  \begin{equation}
    \abs{a+b\tau} < \psi q^2
  \end{equation}
\end{subequations}
with $\psi=4$ have no integer solutions.

With the previous inequality linear forms in logarithms come into play. The
theory to solve this problem, unfortunately, provides only solutions for
huge~$w$ (and fixed $q$), see Section~\ref{sec:bounds-w}. The words
``unfortunately'' and ``huge'' here mean that it is not possible to test the
remaining finitely many configurations in reasonable time. In order to reduce
the bounds from which on \eqref{eq:Main_Log_Form} does not have any integer
solutions (and thus reducing the calculation time), we use a method due to
Baker and Davenport~\cite{Baker-Davenport:1969} in
Section~\ref{sec:smaller-bound}.

We are left with a bunch of small cases. Some remarks on how to check the lemma
for these values directly can be found in Section~\ref{sec:testing-fixed-w}.

Note that the steps above were presented in reverse ordering (from the
perspective, that we only use results, which were proven earlier in the
article), since this is more the way one has to think when solving such a
problem.


\section{Testing Directly}
\label{sec:testing-fixed-w}


In this section, we collect some remarks on how to directly test whether
Proposition~\ref{pro:find-lattice-point} holds for fixed parameters.
So let us fix $q$, $p$ and $w$. We use the following criterion to find a
lattice point $\lambda\in \Lambda_\tau \setminus \Lambda_{\tau^2}$ inside the
rectangle~$R_w$.

We first establish necessary and sufficient conditions for a complex number
$z=x+iy$ to be contained in $R_w$.

\begin{proposition}\label{pro:test-in-Rw-direct}
  Set
  \begin{align*}
    B_1 &= \frac{q^{w+1}}{4}, \\
    B_2 &= \frac{q^{w+1}}{4} + \frac{q^{(w+1)/2}}{2} \sqrt{\frac{q-1/4}{q+2}}
    \intertext{and}
    B_3 &= \frac{q^{w+1}}{2(4q-1)} - \frac{q^{w/2+1}}{\sqrt{4q-1}},
  \end{align*}
  and let us write $\frac12 \tau^{w+1} = u_w + v_w\tau$. Then $\lambda =
  a\tau+b\tau^2 \in R_w$ if and only if
  \begin{subequations}
    \begin{equation}\label{eq:Rw-bound:axial}
      B_1 < 
      a \left(\frac{u_wp}{2}+v_wq\right)
      + b \left(\frac{u_w}{2}-u_wq+\frac{v_wpq}{2}\right)
      < B_2
    \end{equation}
    and
    \begin{equation}\label{eq:Rw-bound:radial}
      \abs{a u_w + b(u_wp + v_wq)}
      < B_3.
    \end{equation}
  \end{subequations}
\end{proposition}

We can solve this system of inequalities and obtain finitely many
pairs of integers~$(a,b)$. If we find a pair with $q\ndivides a$, then
Proposition~\ref{pro:find-lattice-point} is true for this instance. Thus,
Proposition~\ref{pro:test-in-Rw-direct} leads to a ``searching
algorithm'' to solve Proposition~\ref{pro:find-lattice-point} for a
particular parameter set.

\begin{proof}
  Let $\frac12 \tau^{w+1} = x_w + iy_w$.
  By elementary geometry we know that a point $(x,y)$ that lies between the
  upper and lower length sides of $R_w$ satisfies
  \begin{equation*}
    \left(\frac{\abs{\tau}^{w+1}}{2}\right)^2
    < x x_w+y y_w
    < \frac{\abs{\tau}^{w+1}}{2} \left(\frac{\abs{\tau}^{w+1}}{2}+s\right)\!.
  \end{equation*}
  Since
  \begin{equation*}
    x_w + i y_w = u_w+\frac{v_wp}{2} + iv_w\sqrt{q-1/4}
  \end{equation*}
  and
  since we want to have constraints for integers $a$ and $b$ in $(x,y)$ with
  \begin{equation*}
    x + iy = \lambda = a\tau+b\tau^2
    = \frac{ap}{2}+\frac{b}{2}-bq+i(a+bp)\sqrt{q-1/4},
  \end{equation*}
  we obtain
  \begin{equation*}
    x x_w+y y_w = a\left(\frac{u_wp}{2}+\frac{v_w}{4}\right)
    + b\left(\frac{1}{2}-q\right)\left(u_w+\frac{v_wp}{2}\right)
    + (a + bp)v_w \left(q-\frac{1}{4}\right)\!,
  \end{equation*}
  from which~\eqref{eq:Rw-bound:axial} follows.

  The inequality for $(x,y)$ being in between the sides of $R_w$ parallel to
  $\tau^w$ is given by
  \begin{equation*}
    \abs{xy_w-yx_w}
    < \frac{\abs{\tau}^{w+1}}{2} 
    \left(\frac{\abs{v_0-v_1}\abs{\tau}^w}2-\sqrt{q} \right)
    = B_3 \sqrt{q-1/4}
  \end{equation*}
  and, likewise as above, we have
  \begin{equation*}
    x y_w - y x_w = \left(
      a\frac{v_wp}{2} + bv_w\left(\frac{1}{2}-q\right)
      - (a+bp)\left(u_w+\frac{v_wp}{2}\right)
    \right) \sqrt{q-1/4}.
  \end{equation*}
  The inequality~\eqref{eq:Rw-bound:radial} follows.

  Therefore the lattice point $\lambda\in\Lambda_\tau\cap R_w$ satisfies both
  inequalities stated in the lemma and, the other way round, all such points
  are inside $R_w$.
\end{proof}


\section{Huge Bounds for $w$}
\label{sec:bounds-w}


This section deals with showing that the Inequalities~\eqref{eq:Main_Log_Form}
have no integer solutions. We do this by providing a method to find for a fixed
$q$ all $w$ such that~\eqref{eq:Main_Log_Form} is satisfied. More precisely, we
will give a (rather huge) bound on $w$ such that solutions (if any) are only
possible for smaller values.


However, for a single fixed $q$ we still have (too) many possiblities to test
all $w$, see Lemma~\ref{lem:Matveev-bound} below and the text
afterwards. Therefore we will reduce the upper bound of $w$ by using a variant
of the Baker--Davenport method \cite{Baker-Davenport:1969}
in Section~\ref{sec:smaller-bound}.


We can restrict ourselves to the following setting. We may assume $b>0$ since
otherwise $-a$, $-b$ and $w$ would satisfy \eqref{eq:Main_Log_Form}. Moreover,
we may assume that $\gcd(a,b)=1$. If $a$ and $b$ would have a common divisor
$d$ then with $a$, $b$ and $w$ also $a/d$, $b/d$ and $w$ would
satisfy~\eqref{eq:Main_Log_Form}.


As a first step we want to find a bound for $w$ for a fixed integer
$q\geq 2$.


\begin{proposition}\label{pro:expl-comp-bound-ineq}
  For every $q$ there exists an explicitly computable bound $w_q$ such that
  the Inequalities~\eqref{eq:Main_Log_Form} do not have
  any integer solutions with $w\geq w_q$.
\end{proposition}


The following lemma states the precise conditions when solutions
of~\eqref{eq:Main_Log_Form} are
possible. Proposition~\ref{pro:expl-comp-bound-ineq} is a direct consequence of
this result.


\begin{lemma}\label{lem:Matveev-bound}
  Solutions to~\eqref{eq:Main_Log_Form} exist only if
  \begin{multline}\label{Matveev_Bound}
    7.72\cdot10^{13} \log q \log(\sqrt\psi q) 
    \log\left(4.87 w \frac{\max\{3\pi,2\log q\}}{\log(\sqrt\psi q)}\right)\\
    >-\log\chi +\frac{w-4}2 \log q
  \end{multline}
  with $\chi=9$ and $\psi=4$ holds.
  
  In particular, for 
  \begin{align*}
    w&\geq 8.68\cdot 10^{15} \log q \log\log q &&\text{if $q\geq 13$}
    \intertext{and}
    w&\geq 1.973 \cdot 10^{16} &&\text{if $q\in\set{2, 3, \dots, 12}$}
  \end{align*}  
  the Inequalities~\eqref{eq:Main_Log_Form} do not have any integer solutions.
\end{lemma}


It is easy to see that for fixed $q$ the Inequality~\eqref{Matveev_Bound}
cannot hold if $w$ is large. For instance $q=2$ yields $w\leq w_2=8.596\cdot
10^{15}$ or for $q=42$ we obtain $w\leq w_{42} = 2.747\cdot 10^{16}$.
This is one of the key results used in the proof of
our main result, Theorem~\ref{thm:non-opt-large-w}.


In the following, we denote by $\f{h}{\alpha}$ the \emph{absolute logarithmic
  height}, which is defined as follows. Let $\alpha$ be an algebraic number of
degree $d$ and with minimal polynomial
\begin{equation*}
  a_0 \prod_{i=1}^d \left(X-\alpha_i\right),
\end{equation*}
then
\begin{equation*}
  \f{h}{\alpha}=\frac 1d\Big(\log \abs{a_0}+ 
    \sum_{i=1}^d \max\{0, \log\abs{\alpha_i}\}\Big). 
\end{equation*}


For the proof of Lemma~\ref{lem:Matveev-bound} and in view
of~\eqref{eq:Main_Log_Form} we apply the following result due to
Matveev~\cite{Matveev:2000:explicit-lower-bound}.


\begin{theorem}[Theorem~2.2 with $r=1$ in
  \cite{Matveev:2000:explicit-lower-bound}]\label{Matveev:Th}
  Denote by $\alpha_1$, \dots, $\alpha_n$ algebraic numbers, not $0$ nor $1$,
  by $\log\alpha_1$, \dots, $\log\alpha_n$ determinations of their logarithms,
  by~$D$ the degree over $\Q$ of the number field $K =
  \Q(\alpha_1,\ldots,\alpha_n)$, and by $b_1$, \dots, $b_n$ rational
  integers. Furthermore let $\kappa=1$ if $K$ is real and $\kappa=2$
  otherwise. For all integers~$j$ with $1\leq j\leq n$ choose
  \begin{equation*}
    A_j\geq \max\set{D \f{h}{\alpha_j}, \abs{\log\alpha_j}, 0.16},
  \end{equation*}
  and set
  \begin{equation*}
    B=\max \set{1} \cup \set*{\abs{b_j} A_j /A_n}{1\leq j \leq n }.
  \end{equation*}
  Assume that
  \begin{equation*}
    b_1\log \alpha_1+\cdots+b_n \log \alpha_n\not=0.
  \end{equation*}
  Then
  \begin{equation*}
    \log \abs{b_1\log \alpha_1+\cdots+b_n \log \alpha_n}
    \geq -C(n,\kappa)\max\set{1,n/6} C_0 W_0 D^2 \Omega
  \end{equation*}
  with
  \begin{gather*}
    \Omega=A_1\cdots A_n, \\
    C(n,\kappa)= \frac {16}{n!\, \kappa} e^n (2n +1+2 \kappa)(n+2)
    (4(n+1))^{n+1} \left( \frac 12 en\right)^{\kappa}, \\
    C_0= \log\left(e^{4.4n+7}n^{5.5}D^2 \log(eD)\right),
    \quad W_0=\log(1.5eBD \log(eD)).
  \end{gather*}
\end{theorem}


\begin{proof}[Proof of Lemma~\ref{lem:Matveev-bound}]
  We observe that in Matveev's theorem we have $n=3$, $D=4$ and $\kappa=2$
  since we use
  \begin{equation*}
    \alpha_3 = \frac{a+b\tau}{\abs{a+b\tau}},
    \qquad
    \alpha_2 = \frac{\tau}{\abs{\tau}}
    \qquad\text{and}\qquad
    \alpha_3 = i.
  \end{equation*}
  Moreover, we set $b_3=1$, $b_2=-w$ and $b_1=k$.

  Next, let us compute the heights of $\frac{a+b\tau}{\abs{a+b\tau}}$ and
  $\frac{\tau}{\abs{\tau}}$. Let us note that for an imaginary quadratic integer $\alpha$
  the algebraic number $\alpha/|\alpha|$ is a zero of
\begin{equation*}
  \abs{\alpha}^4 (X-\alpha/\abs{\alpha}) (X-\bar\alpha/\abs{\alpha})
  (X+\alpha/\abs{\alpha}) (X+\bar\alpha/\abs{\alpha})
\end{equation*}
and therefore
\begin{equation}\label{eq:heightbound}
  \f{h}{\frac{\alpha}{\abs{\alpha}}}
  = \frac 14\left(\log\abs{a_0} + 4\log\abs{\alpha/\abs{\alpha}}\right)
    \leq \log\abs{\alpha}
\end{equation}
since $\abs{a_0}\leq \abs{\alpha}^4$.
We choose
\begin{align*}
A_3 &= 8 \log (\sqrt\psi q) = 4\log(\psi q^2)\geq 4 \log\abs{a+b\tau},\\
A_2 &= 2 \log q=4\log \abs{\tau},\\
A_1 &= 2\pi=4\abs{\log i}.
\end{align*}

Next we find an upper bound for $k$. Let us note that $\abs{\log \frac{a+b\tau}{\abs{a+b\tau}}}<2\pi$
and from the consideration above we have
\begin{equation*}
  \abs{\log \frac{\tau}{\abs{\tau}}}
  \leq \pi-\arctan(\sqrt{4q-1})\leq \frac{2\pi}3.
\end{equation*}
Therefore, a very crude estimate of inequality~\eqref{eq:Main_Log_Form-a}
yields
\begin{equation*}
  \abs{2\pi+w\frac{2\pi}3-k\frac{\pi}2}<\frac{\pi}2,
\end{equation*}
and, thus, $k\leq \frac{3w}2$. We choose
\begin{equation*}
  B=w \frac{\max\{3\pi,2\log q\}}{A_3}.
\end{equation*}
Before we may apply Theorem~\ref{Matveev:Th} we have to check that our linear
form in logarithms (i.e., the left hand side of~\eqref{eq:Main_Log_Form-a}) is
nonzero. Let us assume for the moment the contrary.
But assuming that the linear form in logarithms is zero expressed in geometric terms is that
$\frac12\tau^{w+1} + \frac12(a\tau+b\tau^2)$ lies on the line segment joining
the points $v_0\tau^w$ and $v_1\tau^w$, thus equals $\frac12\tau^{w+1}$, which
contradicts Lemma~\ref{lem:lower_line}.
In view of inequality~\eqref{eq:Main_Log_Form-a} and Theorem~\ref{Matveev:Th}
we obtain~\eqref{Matveev_Bound}.

We are left to compute the explicit bounds for $w$. Let us assume for the moment that
$\max\{3\pi,2\log q\}=2\log q$, i.e., that $q\geq 112> e^{3 \pi/2}$. Let us note that
under this assumption we have $\log (\sqrt\psi q) < 1.15 \log q$.
By a crude estimate
we deduce that Inequality~\eqref{Matveev_Bound} is not satisfied if
\begin{equation}\label{eq:expl-bound:intermediate}
  1.954 \cdot 10^{14}\log q \log w < w
\end{equation}
holds, unless $w<10^{10}$. Due to a result of Peth\H{o} and
de~Weger~\cite{Pethoe:1986}, namely their Lemma~2.2, an inequality of the form
$A\log x \geq x$ with $A>e^2$ implies the inequality $x<2A\log A$.
Therefore we find an explicit bound for $w$, namely
\begin{equation*}
 w
 \geq 8.68 \cdot 10^{15} \log q \log\log q
 > 3.908 \cdot 10^{14}\log q \log \bigl(1.954 \cdot 10^{14}\log q\bigr),
\end{equation*}
which implies~\eqref{eq:expl-bound:intermediate} and consequently the non-existence of solutions.

By solving inequality~\eqref{Matveev_Bound} for each integer $2\leq q <112$ one can easily show that the explicit bounds stated in the Lemma also
hold for $q<112$.
\end{proof}


\begin{proof}[Proof of Proposition~\ref{pro:expl-comp-bound-ineq}]
  The result follows out of Lemma~\ref{lem:Matveev-bound} since for fixed $q$
  Inequality~\eqref{Matveev_Bound} does not hold if $w$ is sufficiently
  large.
\end{proof}


\section{Reducing the Bounds for $w$}
\label{sec:smaller-bound}


The bounds from the previous section
(Proposition~\ref{pro:expl-comp-bound-ineq}) are too huge in order to test all
remaining configurations in reasonable time. Now our aim is to reduce these
bounds which is done below and works very well: For instance, the bound
$w_2 =8.596\cdot 10^{15}$ is reduced to $\wt w_2=140$. We modify a method due to
Baker and Davenport~\cite{Baker-Davenport:1969} to succeed, see
Lemma~\ref{lem:Baker-Davenport} for details.

The remaining section deals with special cases and the occurrence of some
linear dependence in the linear form in
logarithms. Lemma~\ref{lem:lin-dependence} allows us to test for this linear
dependence and in Lemma~\ref{lem:linear-dependent} we describe how to deal with
this situation. At the end, we deal with two special cases
(Lemma~\ref{lem:lin-dep-a-1}).


\newcommand{\nearest}{\norm}

Let us denote the distance to the nearest integer by $\nearest{\,\cdot\,}$.


\begin{lemma}\label{lem:Baker-Davenport}
  Suppose we have a bound $w_q$ for $w$ (i.e.,
  Inequalities~\eqref{eq:Main_Log_Form} do not have any integer solutions for
  $w\geq w_q$). Fix $a$ and $b$. Let $P/Q$ be a convergent to
  \begin{equation*}
    \epsilon = \frac{2}{i\pi} \log \frac{\tau}{\abs{\tau}}
  \end{equation*}
  with the properties that
  $\nearest{Q\epsilon} = \kappa/w_q$ for some $\kappa<1/4$, but
  $\nearest{Q\delta}>2\kappa$, where
  \begin{equation*}
    \delta = \frac{2}{i\pi} \log \frac{a+b\tau}{\abs{a+b\tau}}.
  \end{equation*}
  Then the Inequalities~\eqref{eq:Main_Log_Form} do not have any
  integer solutions with our fixed $a$ and $b$, and with
  \begin{equation*}
    w \geq \wt w_{a,b} =
    \frac{2}{\log q}\log\left(\frac{2\chi Q}{\kappa\pi}\right)+4,
  \end{equation*}
  where $\chi=9$.
\end{lemma}
Note that $\kappa<1/4$ is formally not needed as an assumption in the
lemma, but $\nearest{Q\delta} > 2\kappa$ implies this condition.

\begin{proof}[Proof of Lemma~\ref{lem:Baker-Davenport}]
  Assume that $w_q > w \geq \wt w_{a,b}$ is satisfied and that we have a solution
  of~\eqref{eq:Main_Log_Form}. We multiply
  Inequality~\eqref{eq:Main_Log_Form-a} by $2Q/\pi$ and use the
  notations of the lemma to obtain
  \begin{equation*}
    \abs{Q\delta - wQ\epsilon + k}
    < \frac{2\chi Q}{\pi}q^{2-w/2}
    \leq \kappa.
  \end{equation*}
  But on the other hand, we have
  \begin{equation*}
    \abs{Q\delta - wQ\epsilon + k}
    \geq \nearest{Q\delta} - w\nearest{Q\epsilon}
    > 2\kappa-w\frac{\kappa}{w_q}
    > \kappa.
  \end{equation*}
  Combining these two inequalities yields a contradiction.
\end{proof}

We want to emphasize that Lemma~\ref{lem:Baker-Davenport} yields bounds only in
case of $\epsilon$ and $\delta$ are linearly independent over $\Q$. If this is
not the case, then the considerations in the remaining section can be used. In particular,
this is the case if $b=0$ or $2a+bp=0$ holds.

The following lemma allows us to test the linear dependence of
$\log\bigl(\frac{\tau}{\abs{\tau}}\bigr)$ and
$\log\bigl(\frac{a+b\tau}{\abs{a+b\tau}}\bigr)$ over $\Q$.

\begin{lemma}\label{lem:lin-dependence}
  Suppose we have integers $a$ and $b$ such that $m = |a+b\tau|^2 <
  \psi^2 q^4$ with $\psi=4$.
  Let us write $q=d_1d_2^2$ and $m=m_1m_2^2$,
  such that $d_2$ and $m_2$ are maximal with respect to
  $\gcd(d_1,d_2)=\gcd(m_1,m_2)=1$.
  Set $m'=m_1/\gcd(m_1,d_1)$. With $\nu_\ell$ being the $\ell$-adic valuation,
  set $\alpha_\ell=\nu_\ell(q)$ for all primes~$\ell\divides q$.
  For odd primes $\ell\divides q$ let $\alpha'_\ell=\alpha_\ell$ if $\ell\divides d_1m'$
  and put $\alpha'_\ell=\alpha_\ell/2$ otherwise.
  If $2\divides q$ we put
  \begin{equation*}
    \alpha'_2=
    \begin{cases}
      \alpha_2/2 &
      \quad\undisp{\text{if $d_1\equiv m'\equiv 1 \!\!\mod 4$,}}\\
      \alpha_2 &
      \quad\undisp{\text{if $2\divides d_1m'$ and
          if $d_1 \equiv 1 \tmod 4$ or $m' \equiv 1 \tmod 4$,}}\\
      \alpha_2 &
      \quad\undisp{\text{if $2\ndivides d_1m'$ and
          if $d_1 \equiv 3 \tmod 4$ or $m' \equiv 3 \tmod 4$,}}\\
      2\alpha_2 &
      \quad\undisp{\text{if $2\divides d_1m'$ and
          if $d_1 \equiv 3 \tmod 4$ or $m' \equiv 3 \tmod 4$.}}
    \end{cases}
  \end{equation*}
  Let $N$ be the greatest common divisor of all $\alpha'_\ell$ with
  primes $\ell\divides q$.

  Then $\log\bigl(\frac{\tau}{\abs{\tau}}\bigr)$ and
  $\log\bigl(\frac{a+b\tau}{\abs{a+b\tau}}\bigr)$
  are linearly dependent over $\Q$ if and only if
  \begin{equation}\label{eq:mult-dependence}
    \left(\frac{\tau}{\abs{\tau}}\right)^\eta
    =\left(\frac{a+b\tau}{\abs{a+b\tau}}\right)^\vartheta
  \end{equation}
  for some positive integer $\vartheta\divides 24N$ and some integer $\eta$
  with $\abs{\eta} < \vartheta\bigl(4+\frac{2 \log \psi}{\log q}\bigr)$.
\end{lemma}

\begin{proof}
  It is immediate that $\log\bigl(\frac{\tau}{\abs{\tau}}\bigr)$ and
  $\log\bigl(\frac{a+b\tau}{\abs{a+b\tau}}\bigr)$ are linearly
  dependent over $\Q$ if and only if \eqref{eq:mult-dependence}
  holds. First, let us consider \eqref{eq:mult-dependence} as an
  equation in the ideal group of the field
  $K=\Q(\tau,\sqrt{d_1},\sqrt{m'})$. 

  We aim to compute the prime ideal factorization of $(\ell)$ for every prime
  $\ell\divides q$. We already know by Lemma~\ref{lem:properties_tau} that
  every such ideal~$(\ell)$ splits in
  $\Q(\tau)$ and is therefore unramified.  Furthermore, by definition
  $d_1$ and $m'$ are coprime; therefore, $(\ell)$ with an odd prime $\ell\divides q$
  is at most in one of the fields $\Q(\sqrt{d_1})$ and $\Q(\sqrt{m'})$
  ramified.  Hence, $(\ell)$ is ramified in $K$ if and
  only if $\ell\divides d_1m'$. Moreover, if the ideal~$(\ell)$ is ramified, then the
  ramification index is exactly two. Altogether we get the following by using
  the fact that $K$ is an Abelian extension of $\Q$.
  \begin{itemize}
  \item If $\ell\ndivides d_1m'$, then $(\ell)=\cJ_\ell \cJbar_\ell$,
    where $\cJ_\ell\divides (\tau)$ is the product of distinct prime
    ideals.
  \item If $\ell\divides d_1 m'$, then $(\ell)=\cJ_\ell^2
    \cJbar_\ell^2$, where $\cJ_\ell\divides (\tau)$ is the product of
    distinct prime ideals.
  \end{itemize}
  
  Let us turn to the case $\ell=2$. Recall that the ideal~$(2)$
  ramifies in the quadratic field $\Q(\sqrt{d})$ if and only if
  $d \equiv 2 \tmod 4$ or $d \equiv 3 \tmod 4$.
  We have to distinguish between several cases.
  \begin{itemize}
  \item First, let us assume that the ideal~$(2)$ neither ramifies in
    $\Q(\sqrt{d_1})$ nor in $\Q(\sqrt{m'})$ (and consequently not in
    $\Q(\sqrt{d_1m'})$ either), i.e., $d_1\equiv m'\equiv 1
    \mod 4$. Then $(2)$ is also unramified in $K$ and we have
    $(2)=\cJ_2 \cJbar_2$, where $\cJ_2$ is the
    product of distinct prime ideals.
  \item There is no situation, when the ideal~$(2)$ ramifies in exactly
    one of the fields $\Q(\sqrt{d_1})$, $\Q(\sqrt{m'})$ and
    $\Q(\sqrt{d_1m'})$.
  \item Next, suppose the ideal~$(2)$ is ramified in exactly
    two of the fields $\Q(\sqrt{d_1})$, $\Q(\sqrt{m'})$ and
    $\Q(\sqrt{d_1m'})$., i.e., one of the two cases
    \begin{itemize}
    \item $2\divides d_1m'$ together with
      $d_1\equiv 1\mod 4$ or $m'\equiv 1 \mod 4$, or
    \item $2\ndivides d_1m'$ together with
      $d_1\equiv 3\mod 4$ or $m'\equiv 3 \mod 4$
    \end{itemize}
    occurs.
    Then
    $(2)=\cJ_2^2 \cJbar_2^2$, where $\cJ_2$ is
    the product of distinct prime ideals.
  \item If the ideal~$(2)$ is ramified in all
    of the fields $\Q(\sqrt{d_1})$, $\Q(\sqrt{m'})$ and
    $\Q(\sqrt{d_1m'})$, i.e., $2\divides d_1m'$ together with 
    $d_1\equiv 3\mod 4$ or $m'\equiv 3\mod 4$, then $(2)=\cJ_2^4 \cJbar_2^4$,
    where $\cJ_2$ is some prime ideal.
  \end{itemize}
  Therefore, by considering the definition of $\alpha'_\ell$, we obtain
  \begin{equation*}
    (\tau) = \prod_{\substack{\ell\divides q\\\text{$\ell$ prime}}}
    \cJ_{\ell}^{2\alpha'_\ell}.
  \end{equation*}
  Note that $\alpha'_\ell=\alpha_\ell/2$ only happens if $\alpha_\ell$
  is even. Furthermore we obtain
  \begin{equation*}
    \left(\frac{\tau}{\abs{\tau}}\right)
    = \prod_{\substack{\ell\divides q\\\text{$\ell$ prime}}}
    \left(\frac{\cJ_\ell}{\cJbar_\ell}\right)^{\alpha'_\ell}\!.
  \end{equation*}
  Therefore $\bigl(\frac{\tau}{\abs{\tau}}\bigr)$ is an $N$th power if
  and only if $N$ divides the greatest common divisor of all
  $\alpha'_\ell$ with primes $\ell\divides q$.

  Let us turn now from the ideal group point of view to the element
  point of view of Equation~\eqref{eq:mult-dependence}.  So far we
  have proved that $\log\bigl(\frac{\tau}{\abs{\tau}}\bigr)$ and
  $\log\bigl(\frac{a+b\tau}{\abs{a+b\tau}}\bigr)$ are linearly
  independent over $\Q$ if and only if there exist integers
  $\vartheta$ and $\eta$ such that
  \begin{equation*}
    \left(\frac{\tau}{\abs{\tau}}\right)^{\eta}
    \left(\frac{a+b\tau}{\abs{a+b\tau}}\right)^{-\vartheta}=1
  \end{equation*}
  and $\vartheta\divides N$. Set $n=\gcd(\eta,\vartheta)$.
  By taking $n$th roots we obtain
  \begin{equation}\label{eq:almost_mult_dependent}
    \left(\frac{\tau}{\abs{\tau}}\right)^{\eta/n}
    \left(\frac{a+b\tau}{\abs{a+b\tau}}\right)^{-\vartheta/n}=\zeta_n,
  \end{equation}
  where $\zeta_n$ is an $n$th root of unity.  The group
  $\Gal(\Q(\zeta_n)/\Q)$ is a subgroup of $\Gal(K/\Q)$, and
  $\Gal(K/\Q)$ is isomorphic to a subgroup of $(\Z/2\Z)^3$. Since
  $24$ is maximal with $\varphi(24)=8$ (where $\varphi$ is Euler's phi
  function), we deduce that $n\divides 24$.
  If we take Equation~\eqref{eq:almost_mult_dependent} to the $24$th
  power, we obtain the first statement of the lemma.
   
  Because of~\eqref{eq:heightbound} we have heights
  \begin{equation*}
    \f{h}{\frac{\tau}{\abs{\tau}}} = \frac{1}{2} \log q
  \end{equation*}
  and
  \begin{equation*}
    \f{h}{\frac{a+b\tau}{\abs{a+b\tau}}}
    \leq \log \abs{a+b\tau}
    = \log \sqrt{m}
    < \log\psi + 2\log q.
  \end{equation*}
  Comparing these heights on the left and right side of
  \eqref{eq:mult-dependence} we obtain
  \begin{equation*}
    \frac{\eta}{2} \log q
    = \eta \f{h}{\frac{\tau}{\abs{\tau}}}
    = \vartheta \f{h}{\frac{a+b\tau}{\abs{a+b\tau}}}
    < \vartheta (\log\psi + 2\log q).
  \end{equation*}
\end{proof}


\begin{lemma}\label{lem:linear-dependent}
  Suppose we have a bound $w_q$ for $w$ (i.e.,
  Inequalities~\eqref{eq:Main_Log_Form} do not have any integer
  solutions for $w\geq w_q$), and suppose that we have for fixed $a$
  and $b$ a linear dependence of the form
  \begin{equation*}
    \eta \log\left(\frac{\tau}{\abs{\tau}}\right)
    = \vartheta\f{\log}{\frac{a+b\tau}{\abs{a+b\tau}}}\!,
  \end{equation*}
  such that $\vartheta>0$.
  \begin{enumerate}
  \item Let $P/Q$ be an expanded fraction of a convergent
    (i.e., $P/Q=P'/Q'$ for a convergent $P'/Q'$)
    to $\epsilon=\frac{2}{\pi}
    \log\bigl(\frac{\tau}{\abs{\tau}}\bigr)$ with the following
    properties. Suppose $W_Q=(Q+\eta)/\vartheta$ is largest
    possible with $W_Q<w_q$ (i.e., $Q<w_q\vartheta-\eta$) and
    $W_Q \in \Z$ such that
    \begin{equation}\label{eq:continued_fraction:approx-dependent}
      Q\abs{\epsilon-P/Q} <\frac{2\chi\vartheta}{\pi}q^{2-W_Q/2}
    \end{equation}
    with $\chi=9$ holds.
    If no such fraction $P/Q$ exists, then set $W_Q=-\infty$.
  \item Let $W$ be the smallest positive integer such that the inequality
  \begin{equation}\label{eq:continued_fraction_valid-II}
    \frac{2\chi\vartheta}{\pi}q^{2-W/2}
    \leq \frac{1}{2(W\vartheta-\eta)}
  \end{equation}
  with $\chi=9$ holds.
  \end{enumerate}
  Then the
  Inequalities~\eqref{eq:Main_Log_Form} do not have any integer
  solutions with our fixed $a$ and $b$, and with
  \begin{equation*}
    w \geq \wt w_{a,b} = \max\{W_Q+1,W\}.
\end{equation*}

\end{lemma}


\begin{proof}
  The assumption on the linear dependence yields an inequality of the form
  \begin{equation}\label{eq:Dioph_two_term-II}
    \abs{(w\vartheta-\eta)\log\left(\frac{\tau}{\abs{\tau}}\right)
      -\vartheta k\frac{i\pi}2}
    < \chi\vartheta q^{2-w/2}
  \end{equation}
  or with the notation of Lemma~\ref{lem:Baker-Davenport},
  \begin{equation}\label{eq:Dioph_two_term-IIb}
    \abs{\epsilon-\frac {\vartheta k}{w\vartheta-\eta}}
    < \frac{2\chi\vartheta}{(w\vartheta-\eta)\pi}q^{2-w/2}.
  \end{equation}
  Note that due to a well-known Theorem of Legendre we have the following: If
  \begin{equation*}
    \frac{2\chi\vartheta}{(w\vartheta-\eta)\pi}q^{2-w/2}
    \leq \frac{1}{2(w\vartheta-\eta)^2},
  \end{equation*}
  which is true for large enough~$w$, then $(\vartheta
  k)/(w\vartheta-\eta)=P'/Q'$, where $P'/Q'$ is a convergent to
  $\epsilon$. Since $P'$ and $Q'$ are coprime, we have $Q'\divides
  w\vartheta-\eta$, so $w\vartheta-\eta=Q$ for some multiple $Q$ of
  $Q'$.
\end{proof}


\begin{lemma}\label{lem:lin-dep-a-1}
  Suppose we have a bound $w_q$ for $w$ (i.e.,
  Inequalities~\eqref{eq:Main_Log_Form} do not have any integer solutions for
  $w\geq w_q$) and suppose that $a=1$ and either $b=0$ or $b=-2p$.
  Let $P/Q$ be an expanded fraction of a convergent
  (i.e., $P/Q=P'/Q'$ for a convergent $P'/Q'$)
 to $\epsilon=\frac{2}{\pi}
  \log\bigl(\frac{\tau}{\abs{\tau}}\bigr)$ with the following
  properties. Suppose $W_Q=Q$ is largest
  possible with $W_Q<w_q$ (i.e., $Q<w_q$) such that
  \begin{equation*}
   Q\abs{\epsilon-P/Q} < \frac{2\chi}{\pi}q^{2-W_Q/2}
  \end{equation*}
  with $\chi=9$ holds.
  If no such convergent exists, then set $W_Q=-\infty$.

  Then the Inequalities~\eqref{eq:Main_Log_Form} do not have any integer
  solutions with our fixed $a$ and $b$ as above, and with
  \begin{equation*}
    w \geq \wt w_{a,b} = \max\{W_Q+1,25\}.
  \end{equation*}
\end{lemma}


Note that the bound~$25$ is sharp for $q=2$, but for $q>2$ a better
bound could be chosen.


\begin{proof}[Proof of Lemma~\ref{lem:lin-dep-a-1}]
  In both cases $\log \bigl(\frac{a+b\tau}{\abs{a+b\tau}}\bigr)$ is an
  integral multiple of $\frac{i\pi}{2}$. Therefore, we consider the inequality
  \begin{equation}\label{eq:Dioph_two_term}
    \abs{w\log\left(\frac{\tau}{\abs{\tau}}\right)-k\frac{i\pi}2}
    < \chi q^{2-w/2}
  \end{equation}
  which is similar to \eqref{eq:Dioph_two_term-II}.
    In the same way as before
  (proof of Lemma~\ref{lem:linear-dependent}) and with the notation of Lemma
  \ref{lem:Baker-Davenport} we obtain
  \begin{equation*}
    \abs{\epsilon-\frac kw} < \frac{2\chi}{w\pi}q^{2-w/2}
  \end{equation*}
  and that $\frac{k}{w}$ equals a convergent $P'/Q'$ to $\epsilon$ if
  \begin{equation}\label{eq:continued_fraction_valid}
    \frac{2\chi}{w\pi}q^{(-w+4)/2}\leq \frac{1}{2w^2}
  \end{equation}
  which is true for large $w$, in particular for all $w\geq 25$.
  Therefore a solution $w\geq 25$ to Inequalities~\eqref{eq:Main_Log_Form}
  corresponds to a fraction $P/Q=P'/Q'$ such that $Q=w$.
\end{proof}


In order to get a reduced bound $\wt w_q$, we look at all possible
combinations of $a$ and $b$ and calculate a bound $\wt w_{a,b}$ by the
lemmata and considerations above. The bound $\wt w_q$ is the maximum
of all these bounds.


\section{Geometry of Numbers}
\label{sec:Geometry_of_Numbers}

\newcommand{\WL}{\mathrm{WL}}
\newcommand{\NH}{\mathrm{NH}}
\newcommand{\RWL}{\wt R_\WL}
\newcommand{\RNH}{\wt R_\NH}
\newcommand{\phiWL}{} 
\newcommand{\phiNH}{} 


The theory of the geometry of numbers is used to show the existence of a
lattice point in the rectangle~$R_w$ with the desired properties (i.e., out of
the lattice $\Lambda_\tau$ but not in $\Lambda_{\tau^2}$). We use two other
rectangles inside~$R_w$, one which is wide but low (called $\RWL$) and
one which is narrow but high (called $\RNH$). Minkowski's lattice point
theorem (Theorem~\ref{thm:minkowski}) gives us the existence of a lattice point
inside each of these two rectangles (see Lemma~\ref{lem:exist-lambda}), and we
are able to construct our desired lattice point out of it
(Lemma~\ref{lem:conditional_existence}), provided that the two found points are
linearly independent. This is guaranteed if the intersection of the two
mentioned rectangles with $\Lambda_{\tau/2}$ only contains $\tau^{w+1}/2$, which follows from the
Inequalities~\eqref{eq:Main_Log_Form} by Lemma~\ref{lem:empty-intersection}. So
much for a short overview on this section; let us begin.


We use the lattices~$\Lambda_\tau$ and~$\Lambda_{\tau/2}$, which were defined in
Section~\ref{sec:lattices}. Throughout this section, we further use the
rectangle $\RNH$ with vertices
\begin{equation*}
  \frac{\tau^{w+1}}2
  \pm i \frac{\tau^{w+1}}{\abs{\tau}^{w+1}}
  q^2 \sqrt{q+2}
\end{equation*}
and
\begin{equation*}
  \frac{\tau^{w+1}}2
  + \frac{\tau^{w+1}}{\abs{\tau}^{w+1}}
  \frac{s \phiNH}{4q}
  \pm i \frac{\tau^{w+1}}{\abs{\tau}^{w+1}}
  q^2 \sqrt{q+2}
\end{equation*}
with $s=\sqrt{\frac{q-1/4}{q+2}}$, and the rectangle
$\RWL$ with vertices
\begin{equation*}
  \frac{\tau^{w+1}}2
  \pm i \frac{\tau^{w+1}}{\abs{\tau}^{w+1}}
  \frac{1}{16\phiWL q}
  \frac{q^{(w+1)/2}}{\sqrt{q-\tfrac14}}
\end{equation*}
and
\begin{equation*}
  \frac{\tau^{w+1}}2
  + \frac{\tau^{w+1}}{\abs{\tau}^{w+1}}
  \cdot 4\phiWL q^{(3-w)/2} \bigl(q-\tfrac14\bigr)
  \pm i \frac{\tau^{w+1}}{\abs{\tau}^{w+1}}
  \frac{1}{16\phiWL q}
  \frac{q^{(w+1)/2}}{\sqrt{q-\tfrac14}}.
\end{equation*}
Note that these two rectangles are both contained in (the closure) of $R_w$.
See also Figure~\ref{fig:small-rectangles}.


In the Lemmata~\ref{lem:empty-intersection}, \ref{lem:points-lambda-in-Rw} and~\ref{lem:conditional_existence} we need that (at least) one of the conditions
\begin{flalign}\label{eq:cond-geom-q-w}
\begin{minipage}{0.8\textwidth}
\begin{itemize}\setlength{\itemsep}{0.5ex}
\item $w\geq 8$ and $q\geq 258$,
\item $w\geq 9$ and $q\geq 17$,
\item $w\geq 10$ and $q\geq 7$,
\item $w\geq 11$ and $q\geq 5$,
\item $w\geq 12$ and $q\geq 4$,
\item $w\geq 13$ and $q\geq 3$ or
\item $w\geq 16$ and $q\geq 2$
\end{itemize}
\end{minipage} &&
\end{flalign}
on $q$ and $w$ holds.
These bounds are sharp in Lemma~\ref{lem:points-lambda-in-Rw}.


\begin{lemma}\label{lem:empty-intersection}
  Suppose $q$ and $w$ satisfy Conditions~\eqref{eq:cond-geom-q-w}.
  If Inequalities~\eqref{eq:Main_Log_Form} do not have any integer solutions (for
  given $q$ and $w$), then the only lattice point of $\Lambda_{\tau/2}$ in
  $\RNH \cap \RWL$ is $\tau^{w+1}/2$.
\end{lemma}


Since by construction $\RNH\cap \RWL$ is a rectangle with side lengths
\begin{equation}\label{eq:size-small-rect}
  2
  q^2 \sqrt{q+2}
  \qquad\text{and}\qquad
  4\phiWL q^{(3-w)/2} \bigl(q-\tfrac14\bigr),
\end{equation}
therefore has an area which decreases
with~$w$, it seems very reasonable to assume that the only lattice point
contained in $\RNH\cap \RWL$ is $\frac12 \tau^{w+1}$. In order to prove this
result, we reformulate this geometric problem into a problem from Diophantine
analysis (finding solutions for Inequalities~\eqref{eq:Main_Log_Form}).


\begin{proof}[Proof of Lemma~\ref{lem:empty-intersection}]
  First let us note that the shortest vector of $\Lambda_{\tau/2}$ is
  $\frac{\tau}2$ which has length $\sqrt{q}/2$ (see Lemma~\ref{lem:shortest_vector}). Therefore the angle
  between the lower long side of $R_w$ and $\lambda\in \RNH\cap \wt
  R_{\mathrm{WL}}$ with $\lambda \neq \frac{\tau^{w+1}}2$ is at most
  \begin{equation}\label{eq:upper_bound_angle}
    \arcsin \frac{4\phiWL q^{(3-w)/2} \bigl(q-\tfrac14\bigr)}{\tfrac12\sqrt{q}}
    = \f[big]{\arcsin}{8\phiWL q^{2-w/2} \bigl(1-\tfrac{1}{4q}\bigr)}
  \end{equation}
  in absolute values. Due to the conditions~\eqref{eq:cond-geom-q-w},
  the argument of the arcsine is less than $0.11$ and we have
  $\arcsin(x)< 9x/8$; we obtain an upper bound for that angle.
  On the other hand the angle between the vector $\tau^w$ and $a+b\tau$ is
  \begin{multline*}
   \abs[big]{\arg{a+b\tau}-\arg{\tau^w}}
   = \abs{\log\!\left(\frac{a+b\tau}{|a+b\tau|}\right)
     - \log\!\left(\frac{\tau^w}{|\tau^w|}\right)}\\
   = \abs{\log\!\left(\frac{a+b\tau}{|a+b\tau|}\right)
     - w\log\!\left(\frac{\tau}{|\tau|}\right) + k\frac{i\pi}2}
  \end{multline*}
  for some integer $k$. This together with
  Inequality~\eqref{eq:upper_bound_angle} yields
  Inequality~\eqref{eq:Main_Log_Form-a}.

  Now let us write
  $\lambda=\frac12\tau^{w+1}+\frac12(a\tau+b\tau^2)$. Further,
  by~\ref{eq:size-small-rect} we know that
  \begin{equation*}
    \abs{\tfrac{\tau}{2} (a+b\tau)}^2 
    \leq q^5 \bigl(1+\tfrac{2}{q}\bigr)
    + 16 q^{5-w} \bigl(1-\tfrac{1}{4q}\bigr)^2
    < 4q^5
  \end{equation*}
  provided that $w\geq 4$. Thus
  \begin{equation*}
    \abs{a+b\tau} < \psi q^2
  \end{equation*}
  with $\psi=4$
  and $a$ and $b$ are bounded in terms of $q$. All together we obtain
  Inequalities~\eqref{eq:Main_Log_Form}.
\end{proof}


In order to find at least one point inside each of the rectangles $\RNH$ and
$\RWL$ we use Minkowski's lattice point theorem (for example, see Theorem~II in
Chapter~III of Cassels~\cite{Cassels:1959:GN}).


\begin{theorem}[Minkowski's lattice point theorem]\label{thm:minkowski}
  Let $S\subset \R^n$ be a compact point set of volume $V$ which is symmetric
  about the origin and convex. Let $\Lambda$ be any $n$-dimensional lattice
  with lattice constant $\f{d}{\Lambda}$. If $V \geq 2^n \f{d}{\Lambda}$ then
  there exists a pair of points $\pm\lambda\in\Lambda\cap S$, with $\lambda
  \neq 0$.
\end{theorem}


\begin{figure}\centering
    \def\lw{12}
  \def\lh{3}
  \def\ax{1.25}
  \def\ay{0.25}
  \def\bx{0.25}
  \def\by{1}

  \begin{tikzpicture}
    [digit/.style={circle, fill=black, inner sep=2pt},]

    \begin{scope}[
    spy using outlines={circle, magnification=3, size=6cm, connect spies}]

    \begin{scope}[scale=0.4, every node/.style={scale=0.33333}, 
      normal/.style={scale=3}]
      \coordinate (v1) at (\lw,0); 
      \coordinate (v0) at (-\lw,0);
      \coordinate (a) at (\ax,\ay);
      \coordinate (b) at (\bx,\by);

      {\draw [fill=lightgray!50] 
        ($(v1)+(-2,0)$) rectangle ($(v0)+(2,\lh)$);
        \node at ($(v0)+(2,\lh)$) [scale=3, anchor=north west] {$R_w$};
      }
      
      {\fill [gray] (-0.5,0) rectangle (0.5,0.5);}

      \draw [thick] (v1)--(v0);
      \draw [dashed, thick] (v1)--($(v1)+(2,-2)$);
      \draw [dashed, thick] (v0)--($(v0)+(-2,-2)$);
      \node [scale=3] at (-10,-1) {digit set $\cD$};

      {\draw (-0.5,0) rectangle (0.5,1.5);
        \node at (0.5,1.5) [anchor=west] {$\RNH$};}
      {\draw (-2,0) rectangle (2,0.5);
        \node at (-2,0.5) [anchor=south west] {$\RWL$};}
      {\node (half) [digit] at (0,0) {};}

      \clip ($(v0)+(-2,-2)$) rectangle ($(v1)+(2,\lh+0.6)$);

      \foreach \c in {-4,-2,...,4}{
        \draw [black!50, dotted] ($-10*(a)+\c*(b)$)--($10*(a)+\c*(b)$);}
      \foreach \c in {-10,-8,...,8}{
        \draw [black!50, dotted] ($-5*(b)+\c*(a)$)--($7*(b)+\c*(a)$);}

      \foreach \a in {-10,-8,...,10}{
        \foreach \b in {-4,-2,...,6}{ 
          \node[draw,circle,inner sep=1pt,fill]
          at (\ax*\a+\bx*\b,\ay*\a+\by*\b) {};
        }}

      {\node (l1) [digit, black] at (a)
      [label={[black, xshift=-3]right:$\lambda_\WL$}] {};
      \draw [->, black, >=latex] (half)--(l1);}
      {\node (l2) [digit, black] at (b)
      [label={[black]left:$\lambda_\NH$}] {};
      \draw [->, black, >=latex] (half)--(l2);}

    \end{scope}

    \node (zoom) at (3.5,-3.2) {};
    \fill[white] (zoom) circle (3cm);
    {\spy [overlay, black!50] on (0,0.17) in node at (zoom);}
    
    \end{scope}

  \end{tikzpicture}

  \caption{The rectangles~$\RNH$ and $\RWL$, and
    points~$\lambda_\NH$ and $\lambda_\WL$.}
\label{fig:small-rectangles}
\end{figure}

\begin{lemma}\label{lem:exist-lambda}
  For $j\in\set{\NH, \WL}$ there exists a lattice point
  $\lambda_j \in \Lambda_{\tau/2}$ in the rectangle~$\wt R_j$.
\end{lemma}

The situation of this lemma is show in Figure~\ref{fig:small-rectangles}.


\begin{proof}
  First, note that the lattice $\Lambda_{\tau/2}$ has lattice constant
  \begin{equation*}
    \f{d}{\Lambda_{\tau/2}} = \tfrac14 \f{d}{\Lambda_\tau}
    = \tfrac14 q\sqrt{q-\tfrac14}.
  \end{equation*}
  Let us mirror the rectangle $\RNH$ on the line joining the points
  $v_0\tau^w$ and $v_1\tau^w$, and consider the rectangle $\RNH$ joint
  with the mirrored rectangle. We obtain a compact, symmetric around
  $\frac12 \tau^{w+1}$, convex set (a rectangle) of volume
  \begin{equation*}
    2\frac{s \phiNH}{4q}
    \cdot
    2
    q^2 \sqrt{q+2}  
    = q\sqrt{q-\tfrac14} = 4 \f{d}{\Lambda_{\tau/2}}.
  \end{equation*}
  Now Minkowski's lattice point theorem yields a $\lambda_\NH \in \RNH\cap
  \Lambda_{\tau/2}$.

  Let us construct $\lambda_\WL$ similarly: Again we mirror the rectangle
  $\RWL$ on the line joining the points $v_0\tau^w$ and $v_1\tau^w$ and
  consider the rectangle $\RWL$ joint with the mirrored rectangle. We obtain a
  compact, symmetric around $\frac12 \tau^{w+1}$, convex set again of volume
  \begin{equation*}
    2 \cdot 4 \phiWL q^{(3-w)/2} \bigl(q-\tfrac14\bigr)
    \cdot
    2 \frac{1}{8\phiWL q}
    \frac{q^{(w+1)/2}}{\sqrt{q-\tfrac14}}
    = q\sqrt{q-\tfrac14} = 4 \f{d}{\Lambda_{\tau/2}}.
  \end{equation*}
  Minkowski's lattice point theorem yields a $\lambda_\WL \in \RWL\cap \Lambda_{\tau/2}$.
\end{proof}


From now on we assume that we have $\lambda_\NH$ and $\lambda_\WL$ as in
Lemma~\ref{lem:exist-lambda}. The following result is needed in the proof of
Lemma~\ref{lem:conditional_existence}. 


\begin{lemma}\label{lem:points-lambda-in-Rw}
  Suppose $q$ and $w$ satisfy Conditions~\eqref{eq:cond-geom-q-w}.
  Then all lattice points of the form
  \begin{equation*}
    \frac{\tau^{w+1}}2
    + a\left( \lambda_\NH-\frac{\tau^{w+1}}2\right)
    + b\left( \lambda_\WL-\frac{\tau^{w+1}}2\right)
  \end{equation*}
  with non-negative integers $a$ and $b$ at most $2q$
  are contained in the rectangle $R_w$.
\end{lemma}


\begin{proof}
  All given points are contained in $R_w$ if the two
  inequalities
  \begin{equation*}
    2q\cdot
    q^2 \sqrt{q+2}
    + 2q\cdot
    \frac{1}{16\phiWL q}
    \frac{q^{(w+1)/2}}{\sqrt{q-\tfrac14}}
    < \frac 12  \sqrt{\frac{q^{w+1}}{4q-1}} -\sqrt{q}
  \end{equation*}
  and
  \begin{equation*}
    2q\cdot
    \frac{s \phiNH}{4q}
    + 2q\cdot
    4\phiWL q^{(3-w)/2} \bigl(q-\tfrac14\bigr)
    < s
  \end{equation*}
  with $s=\sqrt{\frac{q-1/4}{q+2}}$ are satisfied. This is the case for the given conditions.
\end{proof}


With the construction above (and the assumptions of
Lemma~\ref{lem:points-lambda-in-Rw}) we are in a position to prove the
following Lemma.


\begin{lemma}\label{lem:conditional_existence}
  Suppose $q$ and $w$ satisfy Conditions~\eqref{eq:cond-geom-q-w}.
  If
  \begin{equation*}
    \RNH\cap \RWL \cap \Lambda_{\tau/2}
    =\set{\tfrac12 \tau^{w+1}},
\end{equation*}
then there exists a $\lambda \in R_w$ with $\lambda
  \in \Lambda_\tau \setminus \Lambda_{\tau^2}$.
\end{lemma}


\begin{proof}
  First we show that the lattice points $\mu_\NH=\lambda_\NH- \frac12\tau^{w+1}$
  and $\mu_\WL=\lambda_\WL- \frac12\tau^{w+1}$ are linearly independent. Let us
  shift the origin to $\frac12 \tau^{w+1}$ and let us rotate the coordinate
  system such that the long ``lower side'' of the rectangle $R_w$, which
  contains the origin, is parallel to the real axis. In this new coordinate
  system, we write $\hat\lambda_\NH$ and $\hat\lambda_\WL$ for $\lambda_\NH$
  and $\lambda_\WL$ respectively. We have $\abs[big]{\re[big]{\hat\lambda_\NH}}
  < \abs[big]{\re[big]{\hat\lambda_\WL}}$ and $0 < \im[big]{\hat\lambda_\WL} <
  \im[big]{\hat\lambda_\NH}$, i.e., $\lambda_\NH$ and $\lambda_\WL$ are not
  colinear and therefore $\mu_\NH$ and $\mu_\WL$ are linearly independent.

  Since we know now that $\mu_\NH$ and $\mu_\WL$ are linearly independent, there
  exists a basis $\nu_1$, $\nu_2$ of $\Lambda_{\tau/2}$ such that
  \begin{equation*}
    \mu_\NH=\alpha_{11} \nu_1
    \qquad\text{and}\qquad
    \mu_\WL = \alpha_{21} \nu_1+\alpha_{22} \nu_2
  \end{equation*}
  with $0<\alpha_{11}$ and $0\leq \alpha_{21}<\alpha_{22}$. In our case this
  means there exists a basis $\nu_1,\nu_2$ for $\Lambda_{\tau/2}$ such that
  $\nu_1$ and $\nu_2$ are contained in the parallelogram with vertices $0$,
  $\mu_\NH$, $\mu_\WL$ and $\mu_\NH+\mu_\WL$. Moreover, by the assumptions of the lemma
  and Lemma~\ref{lem:points-lambda-in-Rw} all lattice points of the form
  $\lambda=\frac{\tau^{w+1}}2+a\nu_1+b\nu_2\in \Lambda_{\tau/2}$ with $0\leq a\leq
  2q$ and $0\leq b \leq 2q$ are contained in the rectangle $R_w$.

  Now let us write
  \begin{equation*}
    \nu_1=\beta_{11} \tau+\beta_{12}\tau^2
    \qquad\text{and}\qquad
    \nu_2= \beta_{21} \tau+\beta_{22}\tau^2.
  \end{equation*}
  Since $\nu_1$ and $\nu_2$ as well as $\tau$ and $\tau^2$ are
  bases to the same lattice $\Lambda_\tau$ we
  conclude that $\begin{pmatrix} \beta_{11} & \beta_{12} \\ \beta_{21} & \beta_{22} \end{pmatrix} \in
  \mathrm{SL}_2(\Z)$.
  Our aim is to show that there exist non-negative integers $a$ and $b$
  at most $2q$
  such that $\lambda\in\Lambda_\tau$ but $\lambda\not\in\Lambda_{\tau^2}$.
  Setting
  \begin{equation*}
    \tfrac12 \tau^{w+1} = u_w\frac{\tau}2+v_w\frac{\tau^2}2,
  \end{equation*}
  it suffices to prove that
  \begin{equation*}
    \left(u_w+a\beta_{11}+b\beta_{12}\right)\frac{\tau}2 +
    \left(v_w+a\beta_{21}+b\beta_{22}\right)\frac{\tau^2}2
    = \gamma_1\frac{\tau}2+\gamma_2\frac{\tau^2}2
  \end{equation*}
  for some $\gamma_1\equiv 2 \mod 2q$ and $\gamma_2\equiv 0 \mod 2q$
  has a solution. But a
  solution can be found from a solution to the linear system
  \begin{align*}
    u_w+a\beta_{11}+b\beta_{12}&=2\\
    v_w+a\beta_{21}+b\beta_{22}&=0
  \end{align*}
  modulo $2q$. Such a solution certainly exists since $\begin{pmatrix}
      \beta_{11} & \beta_{12} \\ \beta_{21} & \beta_{22} \end{pmatrix} \in
  \mathrm{SL}_2(\Z)$.
\end{proof} 


With the previous results, we are ready to prove the following proposition.


\begin{proposition}\label{pro:expl-comp-bound-Rw}
  For every $q$ there exists an explicitly computable bound $w_q$ such that
  for each $w\geq w_q$ there exists a lattice point
  \begin{equation*}
    a \tau + b (q-p\tau)
  \end{equation*}
  with $a\in\Z$, $q\ndivides a$ and with $b\in\Z$ in the rectangle~$R_w$.
\end{proposition}


\begin{proof}
  Proposition~\ref{pro:expl-comp-bound-ineq} states a similar result for the
  non-existence of solutions of Inequalities~\eqref{eq:Main_Log_Form}. By
  Lemma~\ref{lem:empty-intersection} this translates to having the
  single intersection point $\frac12 \tau^{w+1}$
  in the rectangles $\RNH$ and $\RWL$. This condition is then
  used in Lemma~\ref{lem:conditional_existence} to find a lattice point
  inside~$R_w$ as desired.
\end{proof}


\section{An Algorithm to Test for fixed $q$}
\label{sec:algorithm}


In order to prove Proposition~\ref{pro:find-lattice-point} for given $q$ and
$p$, but all $w$---this means showing the existence of a lattice point
in each rectangle~$R_w$ with the stated properties---we apply the
following algorithm\footnote{See
  \url{http://www.danielkrenn.at/koblitz2-non-optimal} for the code.}


\begin{algo}\label{alg:test-fixed-q}
  We fix $q$ and fix a choice of $p\in\set{-1,1}$ as input. This
  algorithm returns a list of values for $w$ for which no lattice point in
  $R_w$ exists. We proceed as follows.
  \begin{enumerate}[leftmargin=2em]
  \item\label{enu:algo:Matveev}
    Compute an upper bound $w_q$ for possible solutions (with $w<w_q$)
    using Matveev's Theorem~\ref{Matveev:Th} and in particular
    Inequality~\eqref{Matveev_Bound}.
  \item\label{enu:algo:reduce}
    Reduce the bound $w_q$ by the Baker--Davenport method
    (Section~\ref{sec:smaller-bound}).
  \begin{enumerate}[leftmargin=1.5em]
  \item\label{enu:algo:precomp}
    Compute sufficiently many\footnote{``sufficiently many'' means
      that in step~\ref{enu:algo:non-dependent} of the algorithm a
      convergent can be found for all (non-dependent) situations $a$
      and $b$.}  (consecutive) convergents $P/Q$ to
    $\epsilon=\frac{2}{i\pi} \log \frac{\tau}{\abs{\tau}}$ and save
    them in a list~$\mathcal L$. Precalculate and save $\kappa$ with
    $\nearest{Q\epsilon}=\kappa/w_q$ as well.
  \item\label{enu:algo:special} Use Lemma~\ref{lem:lin-dep-a-1} to deal
    with the case $2a+bp=0$ (i.e., $a=1$, $b=-2p$) and
    compute the new bounds~$\wt w_{a,b}$.
  \item For all integers $a$, $b$ with
    $\abs{a+b\tau} < \psi q^2$, $\psi=4$ and with $b>0$, coprime $a$, $b$ and
    excluding the situations from step~\ref{enu:algo:special} do the following:
  \begin{enumerate}[leftmargin=1.5em]
  \item\label{enu:algo:eta-theta} Find $\eta$ and $\vartheta$ such that
    \begin{equation*}
      \left(\frac{\tau}{\abs{\tau}}\right)^\eta
      =\left(\frac{a+b\tau}{\abs{a+b\tau}}\right)^\vartheta
    \end{equation*}
    using Lemma~\ref{lem:lin-dependence}.
  \item\label{enu:algo:non-dependent}
    If such $\eta$ and $\vartheta$ do not exist in
    step~\ref{enu:algo:eta-theta}, find the convergents~$P/Q$ in
    $\mathcal L$ with smallest $Q$ that satisfies $\kappa<\frac{1}{4}$ and
    $\nearest{Q\delta}>2\kappa$ with
    $\delta = \frac{2}{i\pi} \log \frac{a+b\tau}{\abs{a+b\tau}}$.
    Compute the new bound~$\wt w_{a,b}$ due to Lemma~\ref{lem:Baker-Davenport}.
  \item\label{enu:algo:dependent}
    If such $\eta$ and $\vartheta$ in
    step~\ref{enu:algo:eta-theta} exist, find an expanded fraction $P/Q$
    of the convergents in
    $\mathcal L$ with largest~$Q$ that satisfies $W_Q=(Q+\eta)/\vartheta\leq w_q$ and \eqref{eq:continued_fraction_valid-II}.
    Compute the new bound~$\wt w_{a,b}$ due to Lemma~\ref{lem:linear-dependent}.
  \end{enumerate}
  \item Calculate $\wt w_q$ as the maximum of all $\wt w_{a,b}$.
  \end{enumerate}
\item\label{enu:algo:test-Rw} For all $4 \leq w < \wt w_q$ (with
  exceptions obtained by taking into account the assumptions of
  Proposition~\ref{pro:find-lattice-point}), we verify
  Proposition~\ref{pro:find-lattice-point} directly as described in
  Section~\ref{sec:testing-fixed-w}.
  \end{enumerate}  
\end{algo}


Note that Algorithm~\ref{alg:test-fixed-q} as it is written is not
guaranteed to terminate.\footnote{One might call
  Algorithm~\ref{alg:test-fixed-q} only a ``procedure''.} The reason is
that it might be impossible to find a convergent~$P/Q$ with the
desired properties in step~\ref{enu:algo:non-dependent}. Stopping this
search at some point and not using the reduced bound of
step~\ref{enu:algo:reduce} will make the algorithm terminate for
sure. However, a huge amount of $w$ have to be checked in
step~\ref{enu:algo:test-Rw} then.


\begin{proposition}
  Let $q\geq2$. If Algorithm~\ref{alg:test-fixed-q} terminates, then
  it is correct, i.e., it returns a list of values for $w$ for which
  no lattice point in $R_w$ with the properties stated in
  Proposition~\ref{pro:find-lattice-point} exists.
\end{proposition}


\begin{proof}
  Section~\ref{sec:Geometry_of_Numbers} reduces the problem of finding
  lattice points in $R_w$ to showing that
  Inequalities~\eqref{eq:Main_Log_Form} do not have any integer
  solutions. Step~\ref{enu:algo:Matveev} provides a bound for~$w$; it
  can be computed effectively according to
  Proposition~\ref{pro:expl-comp-bound-ineq}. Step~\ref{enu:algo:reduce}
  reduces this bound. We get a bound for each possible combination of
  $a$ and $b$ (all the different cases are analyzed in
  Section~\ref{sec:smaller-bound}); correctly determining whether we
  have a linear dependence is done via Lemma~\ref{lem:lin-dependence}.

  Taking the maximum of all these bounds yields $\wt w_q$, a new
  bound. In step~\ref{enu:algo:test-Rw} all remaining $w$ are checked
  by a direct search according to
  Proposition~\ref{pro:test-in-Rw-direct}. Since this proposition
  finds a lattice point if and only if one exists, we are able to
  classify all the $w$ and return a list of the exceptional values.
\end{proof}


\begin{proof}[Proof of Proposition~\ref{pro:find-lattice-point}]
  We apply Algorithm~\ref{alg:test-fixed-q}.
\end{proof}


\part{The Part with the Digits}
\label{sec:digits}


\section{Overview}
\label{sec:digits-overview}


In this part of the article, we construct the actual counterexamples to the
minimality of the width\nbd-$w$ non-adjacent forms (see
Section~\ref{sec:overview} for the relevant definitions). This means we have to
find an expansion of a lattice point with a lower number of nonzero digits
than the width\nbd-$w$ non-adjacent form (\wNAF{}) of this point.

We reuse the ideas of Heuberger and
Krenn~\cite{Heuberger-Krenn:2013:wnafs-optimality} for our construction. This
work also tells us that (if it exists) a counterexample using a
(multi-)expansion of weight two can be found. Therefore, we will
try to find
\begin{equation}\label{eq:nonopt-general}
  A\tau^{w-1} + B = C\tau^w + D = E\tau^{2w} + F\tau^w + G,
\end{equation}
where the most left and the most right parts of the equation are valid digit expansions,
i.e., $A$, $B$, $E$, $F$ and $G$ are digits. Moreover, we assume that $D$ is a digit (equal to
$G$), but, in order to get a counterexample to minimality, $C$ is not allowed
to be a digit. However, the point~$C$ is important during the construction of this
counterexample: we will have $C=E\tau^w+F$ and $D=G$, and, more important,
there will be a change $\Delta$ with $\tau C = A+\Delta$ and $D = B -
\Delta\tau^{w-1}$.

Some explicit constructions are given in Section~\ref{sec:non-opt-w}
and Proposition~\ref{pro:nonopt-q4w6}, but most
of the time we will consider a more general situation. There, the existence of
a construction like above relies on Proposition~\ref{pro:find-lattice-point}, which
was proven in the previous part. This lemma gives us the point~$C$. The
change~$\Delta$ is discussed in Section~\ref{sec:small-change-Delta}, and
Section~\ref{sec:lsd} deals with the digits $B$, $D$ and $G$. Everything is
glued together in Sections~\ref{sec:strategy} and~\ref{sec:proof-main}.

We begin with a section which deals with the digit set we
use. Note that this digit set is strongly related to the Voronoi cell defined
in Section~\ref{sec:setup}.


\section{Digit Sets}
\label{sec:digit-sets}


In this section, we make a formal definition of the used digit set. This is
equivalent to the definition stated in Section~\ref{sec:overview} but uses the
Voronoi cell to model the minimal norm property. Afterwards we show that this
choice of digits is unique.


\begin{definition}[Minimal Norm Digit Set]
  \label{def:min-norm-digit-set}

  Let $w$ be an integer with $w\geq2$ and $\cD\subseteq\Ztau$ consist of $0$
  and exactly one representative of each residue class of $\Ztau$ modulo
  $\tau^w$ that is not divisible by $\tau$. If all such representatives
  $\eta\in\cD$ satisfy $\eta \in \tau^w V$, then $\cD$ is called the
  \emph{minimal norm digit set modulo $\tau^w$}.
\end{definition}


The minimal norm digit set above is uniquely determined, see below.


\begin{proposition}\label{pro:digit-set-unique}
  Let $\cD$ be a minimal norm digit set modulo $\tau^w$. Then $\cD$ is uniquely
  determined. In particular, there exists a unique element of minimal norm in
  each residue class modulo $\tau^w$ which is not divisible by~$\tau$.
\end{proposition}


This proposition was proved for $q=2$ in Avanzi, Heuberger and
Prodinger~\cite{Avanzi-Heuberger-Prodinger:2011:redun-expan-i}. The proof there
uses a result of Meier and Staffelbach~\cite{Meier-Staffelbach:1993:effic},
namely their Lemma~2. This lemma and the result for $q=2$ can be generalized in
a straight forward way for arbitrary primes~$q$. We use a different method in
this article, which gives us the result for arbitrary integers~$q$ here.

\begin{proof}[Proof of Proposition~\ref{pro:digit-set-unique}]
  Digits strictly inside the scaled Voronoi cell~$\tau^wV$ are unique, since
  they are closer to $0$ than to any other point of $\tau^w\Ztau$ by the
  definition of the Voronoi cell. In
  Lemma~\ref{lem:lower_line} we have already shown that there are no lattice
  points on the boundary of the scaled Voronoi cell~$\tau^wV$. Therefore no
  non-uniqueness can occur and thus the proposition is proved.
\end{proof}


\section{Non-Optimality for some Values of $w$}
\label{sec:non-opt-w}


We start here with a first family of counterexamples to optimality. We show
the existence of expansions like in~\eqref{eq:nonopt-general}. The following
propositions are devoted to the case when $w=2$, where we give an explicit
construction. Afterwards, we consider the case $w=3$.


\begin{proposition}\label{pro:nonopt-w2-pp}
  Let $q\geq3$ and $p=1$, and set $a=\ceil{q/2}$. Set
  \begin{align*}
    & &
    E &= 1, \\
    A &= (1-a)\tau+a-q, &
    F &= q-a, \\
    B &= 1-\tau, &
    G &= (a-1)\tau+1.
  \end{align*}
  Then
  \begin{equation}\label{eq:nonopt-w2-pp}
    A\tau + B = E\tau^4 + F\tau^2 + G
  \end{equation}
  and both sides of the equation are valid digit expansions, i.e., the $2$-NAF
  is not a minimal digit expansion.
\end{proposition}


\begin{proposition}\label{pro:nonopt-w2-pm}
  Let $q\geq3$ and $p=-1$, and set $a=\ceil{q/2}$. Set
  \begin{align*}
    & &
    E &= 1, \\
    A &= (1-a)\tau-a+q, &
    F &= q-a, \\
    B &= (a-1)\tau-1, &
    G &= -\tau-1.
  \end{align*}
  Then
  \begin{equation}\label{eq:nonopt-w2-pm}
    A\tau + B = E\tau^4 + F\tau^2 + G
  \end{equation}
  and both sides of the equation are valid digit expansions, i.e., the $2$-NAF
  is not a minimal digit expansion.
\end{proposition}


\begin{proof}[Proof of Propositions~\ref{pro:nonopt-w2-pp}
  and~\ref{pro:nonopt-w2-pm}]
  This proof is assisted\footnote{The worksheet can be found at
    \url{http://www.danielkrenn.at/koblitz2-non-optimal}.}
  by SageMath~\cite{SageMath:2016:7.0}.
  Equality
  in~\eqref{eq:nonopt-w2-pp} and~\eqref{eq:nonopt-w2-pm} is easy to verify and
  can be done by a simple symbolic calculation over the ring $\Z[\tau]$.

  We are left with checking that we have valid digit expansions, i.e., that all
  claimed digits are indeed digits. We do this by showing that these quantities
  are closer to $0$ than to any neighbouring lattice point of $\tau^2\Z[\tau]$
  (see also the construction of the digit set via Voronoi cells,
  Section~\ref{sec:digit-sets}). These neighboring lattice points are exactly
  the $\tau^2$-multiples of the points $1$, $\tau$, $\tau-p$, $-1$, $-\tau$ and
  $-\tau+p$. This leads to six inequalities for each digit. Note that we have
  $\abs{a+b\tau}^2 = a^2+b^2q+pab$. We also check that the point $C =
  \tau^2+(q-a) = p\tau-a$ is not a digit for technical reasons, which leads to
  one additional inequality.  Note that we have $\Delta=-pa$ (with the notation
  of Section~\ref{sec:digits-overview}) here.

  However, distinguishing between $p=1$ and $p=-1$ and between $q=2\wt q$
  (even) and $q=2\wt q-1$ (odd), we get $25$ polynomials (each as difference
  of the two sides of an inequality) out of $\Z[\wt q]$. All these
  polynomials have degree at most $3$ and a positive leading coefficient, and
  we can show, by using interval-arithmetic, that all their roots are smaller
  than $2$. This means all polynomials are positive for $\wt q \geq 2$ and
  therefore the inequalities are satisfied.
  
  Since the constant terms of the claimed digits are not divisible by $q$, they
  themselves are not divisible by $\tau$. Therefore, we get valid digit
  expansions, which finishes the proof.  
\end{proof}


\begin{proposition}\label{pro:nonopt-w3-pp}
  Let $q\geq3$ and $p=1$, and set $a=\ceil{q/2}$ and $b=\ceil{q^2/(6q-2)}$. For
  odd $q$ set
  \begin{align*}
    & &
    E &= 1, \\
    A &= (1-q)\tau + (a-1)q-1, &
    F &= (q-a)\tau + q-a, \\
    B &= (q-a)\tau + a-1, &
    G &= -a\tau + a-1,
    \intertext{and for even $q$ set}
    & &
    E &= 1, \\
    A &= (-a-b)\tau+(a-1)q+1, &
    F &= (q-a-1)\tau+q-b, \\
    B &= a\tau-aq+1, &
    G &= -(a-1)\tau+(a-1)q+1.
  \end{align*}
  Then
  \begin{equation*} 
    A\tau^2 + B = E\tau^6 + F\tau^3 + G
  \end{equation*}
  and both sides of the equation are valid digit expansions, i.e., the $3$-NAF
  is not a minimal digit expansion.
\end{proposition}


\begin{proposition}\label{pro:nonopt-w3-pm}
  Let $q\geq3$ and $p=-1$, and set $a=\ceil{q/2}$ and $b=\ceil{q^2/(6q-2)}$. For
  odd $q$ set
  \begin{align*}
    & &
    E &= 1, \\
    A &= (q-1)\tau + (a-1)q-1, &
    F &= (q-a)\tau + q-a, \\
    B &= -a\tau + 1-a, &
    G &= (q-a)\tau + 1-a,
    \intertext{and for even $q$ set}
    & &
    E &= 1, \\
    A &= (a+b)\tau+(a-1)q+1, &
    F &= (q-a-1)\tau-q+b, \\
    B &= -a\tau-aq+1, &
    G &= (a-)\tau+(a-1)q+1.
  \end{align*}
  Then
  \begin{equation*} 
    A\tau^2 + B = E\tau^6 + F\tau^3 + G
  \end{equation*}
  and both sides of the equation are valid digit expansions, i.e., the $3$-NAF
  is not a minimal digit expansion.
\end{proposition}


\begin{proof}[Proof of Propositions~\ref{pro:nonopt-w3-pp}
  and~\ref{pro:nonopt-w3-pm}]
  We use the same machinery as in the proof of
  Propositions~\ref{pro:nonopt-w2-pp} and~\ref{pro:nonopt-w2-pm}.

  If $p=1$ and $q$ is odd, then we take $C = (1-a)\tau - a$ for our technical
  point. We use $\Delta = 1-\tau$. Note that $B=a-1 + (q-a)\tau$ and $D=G=a-1 -
  a\tau$. To verify that $A$, $B$, $E$, $F$ and $G$ are digits, we, again,
  calculate the distance to $0$ and its neighbours in $\tau^3\Z[\tau]$. This
  results in $31$ inequalities, which are of polynomial type. This leaves us
  to check if elements out of $\Z[\wt q]$ with degree at most~$4$ and
  $q=2\wt q-1$ are positive, see the proof above for details. We can affirm
  this (it was done algorithmically).

  For $p=-1$ and odd $q$ we use $C = (1-a)\tau + a$, $\Delta = -\tau-1$, $B=1-a
  - a\tau$, $D=G=1-a + (q-a)\tau$ and proceed in the same manner.

  In the case that $p=1$ and $q$ is even, we have to be more careful because of
  the definition of $b$. We start similar and take $C = -a\tau-b$ for our
  technical point and use $\Delta = q-1$. Moreover, we use $B = 1 +
  a\tau^2$ and $D = G = 1 - (a-1)\tau^2$. The verification of the digits is
  done as above, but we take $b$ into account.
  
  Since $b=\ceil{q^2/(6q-2)}$ is obviously not polynomial, we cannot
  expect that these distance inequalities are polynomials.  To deal
  with the ceil-rounding, we use for the moment $b_1 = q^2/(6q-2)+1$ and
  $C_1=-a\tau-b_1$
  (note, this is not a lattice point) instead of $b$ and $C$ in the
  resulting inequalities; we will correct this later. As this $b_1$ is
  rational, we multiply the inequalities first by $6q-2$ and then
  check whether the resulting polynomials (difference of the two sides
  of the inequality) out of $\Z[\wt q]$ are positive for all $q=2\wt
  q$ as in the proof above. This verification is successful. This
  particularly means, that the point $-a\tau-b_1$ is not in $\tau^3 V$
  (the scaled Voronoi cell containing the digits).

  Next, consider $C_0 = - a \tau - b_0$ with $b_0=q^2/(6q-2)$. We have
  \begin{equation*}
    C_0 = \frac{\tau^3}{2} + ci\tau^3
  \end{equation*}
  with $c = \im{\tau}/(3q-1)$. This means that $C_0$ is on the
  boundary of $\tau^3 V$ (as $\tau^3/2$ is on the boundary).
  The point~$C=-a\tau-b$ is located 
  on the line from $C_0$ to $C_1$. Due to convexity of $\tau^3 V$, this
  lattice point~$C$ lies on the outside of $\tau^3 V$ and, thus, is not
  a digit.

  If we have $p=-1$, still with an even $q$, the proof works similarly, but we
  use $C= -a\tau+b$.
\end{proof}


\section{Existence of a Small Change}
\label{sec:small-change-Delta}


From now on, in contrast to the previous section, where we have given an
explicit construction for a counterexample to minimality of the \wNAF{}, we
start with a different approach. It still builds up on the ideas mentioned in
the introduction of Part~\ref{sec:digits} and will be described fully in
Section~\ref{sec:strategy}. Before we are ready for this alternative
construction, we need a couple of auxiliary results. In this section, we look
at the change $\Delta$ a bit more closely, see
Lemma~\ref{lem:point-in-int-V}, but let us start with the following lemma.


\begin{lemma}\label{lem:parallelogram}
  Let $P$ be the parallelogram with vertices $1$, $\tau-p$, $-1$,
  $-\tau+p$. Then a disc with center $0$ and radius
  $s=\sqrt{\frac{q-1/4}{q+2}}$ fits exactly (i.e., the radius~$s$ is largest
  possible) in the parallelogram~$P$.
\end{lemma}

\begin{figure}
  \centering
  \includegraphics{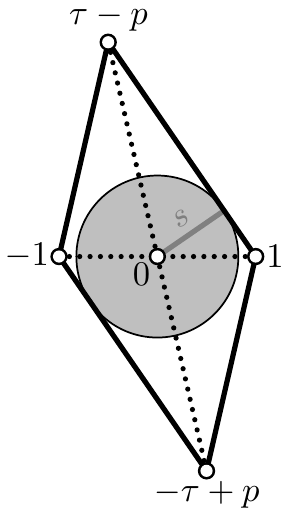}
  \caption{Parallelogram~$P$ of Lemma~\ref{lem:parallelogram} for $q=5$ and
    $p=1$.}
  \label{fig:parallelogram-circle}
\end{figure}

The situation is shown in Figure~\ref{fig:parallelogram-circle}. Note that we
have $s \geq \sqrt{7}/4$.


\begin{proof}
  We can assume $p=1$, since the other situation ($p=-1$) is just
  mirrored. First, we calculate the difference of the areas of the two
  triangles with vertices $\tau-1$, $-\frac12$, $1$ and $\tau-1$, $-\frac12$,
  $0$ and get
  \begin{equation*}
    \tfrac12 \im{\tau} \left(\tfrac32-\tfrac12\right) 
    = \tfrac{1}{2} \sqrt{q-\tfrac14}.
  \end{equation*}
  This area is equal to the area of the triangle with vertices $0$, $1$ and
  $\tau-1$, which is
  \begin{equation*}
    \tfrac12 s \sqrt{\im\tau^2+\tfrac94} = \tfrac12 s \sqrt{q+2}.
  \end{equation*}
  Therefore the desired $s$ follows.
\end{proof}


\begin{lemma}\label{lem:point-in-int-V}
  Let $z$ be in the rectangle~$R_w$. Then there exists a
  \begin{equation*}
    \Delta\in\set{-1,1,\tau-p,-\tau+p},
  \end{equation*}
  such that $z-\Delta$ is in the interior of the scaled Voronoi cell $\tau^wV$.
\end{lemma}


Note that we will have $A=z-\Delta$ and $z=\tau C$ in the construction of our
expansions. Moreover, the point $z$ is the point inside the rectangle~$R_w$
whose existence was shown as the main result of Part~\ref{sec:dioph},
Proposition~\ref{pro:find-lattice-point}.


\begin{proof}[Proof of Lemma~\ref{lem:point-in-int-V}]
  Let $P_z$ be the parallelogram with vertices $z+1$, $z+\tau-p$, $z-1$,
  $z-\tau+p$, i.e., a shifted version of the parallelogram~$P$ of
  Lemma~\ref{lem:parallelogram}, see also
  Figure~\ref{fig:parallelogram-circle}. Since $z$ is in the (open)
  rectangle~$R_w$, its distance to the line~$L$ from~$\tau^w v_0$ to~$\tau^w
  v_1$ is smaller than~$s$ (note that $s$ is the height of the rectangle~$R_w$,
  Section~\ref{sec:setup}). Further, since the rectangle~$R_w$ starts
  $\sqrt{q}$ away from the points~$\tau^w v_0$ and~$\tau^w v_1$ respectively,
  the distance from~$z$ to one of those two points is larger than
  $\sqrt{q}$. Therefore, the line~$L$ cuts the parallelogram~$P_z$ into two
  parts. The two cutting points are on different edges of~$P_z$. This means,
  that there exists a vertex of the parallelogram~$P_z$ on that side of the
  line~$L$, where there is no rectangle~$R_w$. Such a point lies in the
  interior of the Voronoi cell~$\tau^wV$, which can be seen using some
  properties of~$\tau^wV$, including that two neighbouring edges of~$\tau^wV$
  have an obtuse angle at their point of intersection and that a disc with
  center~$0$ and radius $\sqrt{q}$ is contained in $\tau^wV$.
\end{proof}


\section{Change in least significant digit}
\label{sec:lsd}


For our construction of the counterexample, we also have to deal with the
change in the least significant digit, i.e., with the digits $B$, $D$ and $G$
in the expansions~\eqref{eq:nonopt-general}.


\begin{figure}
  \centering
  \includegraphics[width=\linewidth]{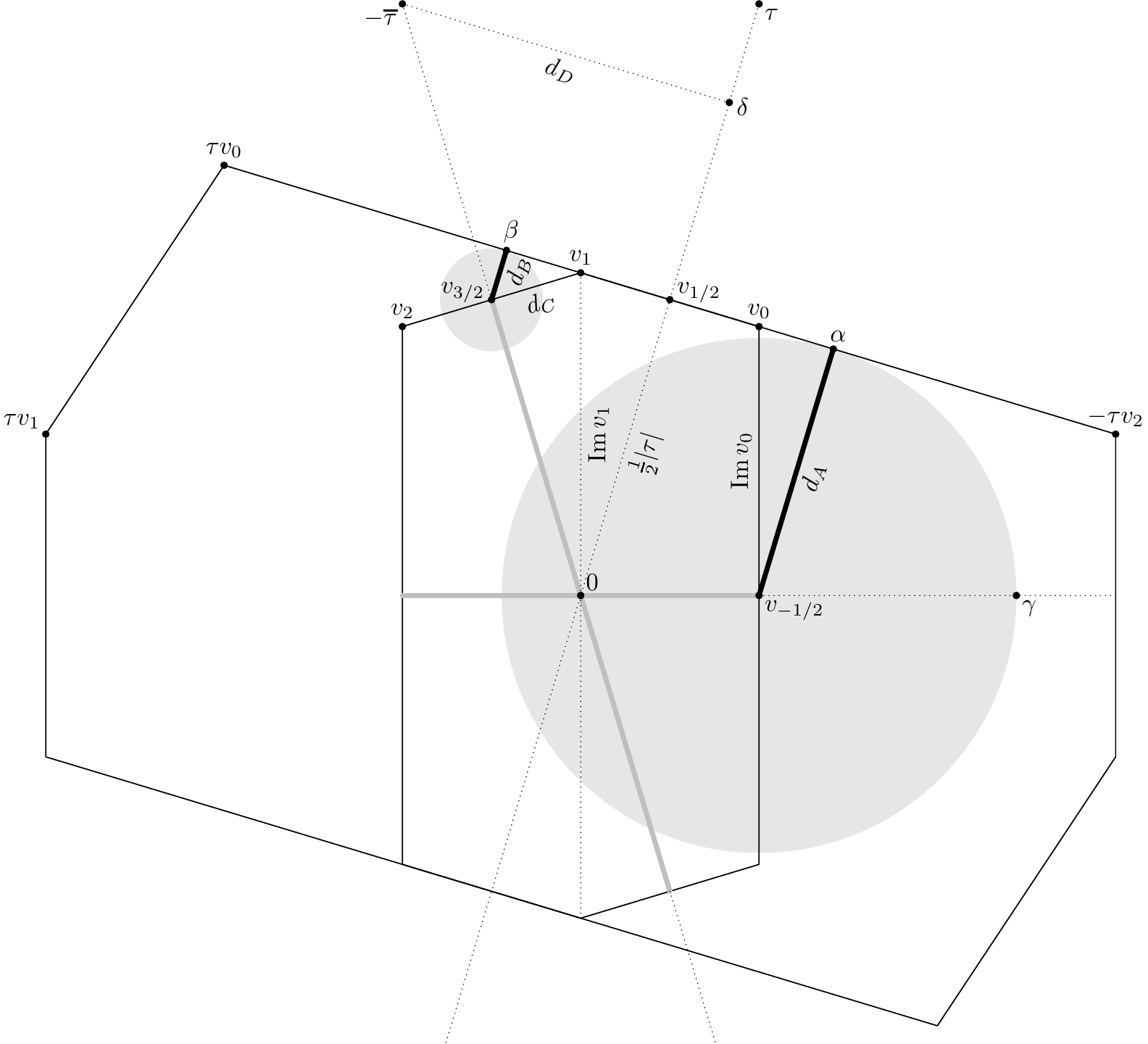}
  \caption{Distances and points used in Lemma~\ref{lem:distances-voronoi-lsd}
    and Proposition~\ref{pro:compensate}.}
  \label{fig:change-lsd}
\end{figure}


\begin{lemma}\label{lem:distances-voronoi-lsd}
  We get
  \begin{equation*}
    d_A = \frac{\sqrt{q}}{2} - \frac{1}{4\sqrt{q}}
  \end{equation*}
  and 
  \begin{equation*}
    d_B = \frac{1}{4\sqrt{q}}
  \end{equation*}
  in Figure~\ref{fig:change-lsd}.
\end{lemma}


\begin{proof}
  The triangle $v_{-1/2} \alpha v_0$ is similar to the triangle $0 v_{1/2}
  v_1$. Therefore
  \begin{equation*}
    d_A = \im*{v_0} \frac{\frac12\abs{\tau}}{\im*{v_1}}
    = \frac{q-\frac12}{2q} \sqrt{q} 
    = \frac{\sqrt{q}}{2} - \frac{1}{4\sqrt{q}}.
  \end{equation*}
  
  To calculate $d_B$ we start as above. The triangle $v_{3/2} \beta v_1$ is
  similar to the triangle $-\conj{\tau} \delta 0$, so $d_B = d_D d_C /
  \abs{-\conj{\tau}}$. We have
  \begin{equation*}
    d_C = \frac12 \sqrt{\left(\im*{v_1}-\im*{v_2}\right)^2 
      + \left(\frac12\right)^2}
    = \frac12 \sqrt{\frac14 + \frac{1}{16(q-\frac14)}}
    = \frac14 \sqrt{1 + \frac{1}{4q-1}}
  \end{equation*}
  by the Pythagorean theorem. The distance~$d_D$ is the projection of
  $-\conj{\tau}$ on the normalized vector with direction $i\tau$. Therefore
  \begin{equation*}
    d_D = \re{\conj{-\conj{\tau}} i \frac{\tau}{\abs{\tau}}}
    = - \frac{1}{\abs{\tau}} \re{i\tau^2}
    = \frac{1}{\abs{\tau}} \im{\tau^2}
    = \frac{1}{\sqrt{q}} \frac12\sqrt{4q-1}
    = \sqrt{1 - \frac{1}{4q}}.
  \end{equation*}
  Now we can calculate $d_B$ as
  \begin{equation*}
    d_B = \frac{1}{\sqrt{q}} \sqrt{1 - \frac{1}{4q}}
    \frac14 \sqrt{1 + \frac{1}{4q-1}}
    = \frac{1}{4\sqrt{q}}.
  \end{equation*}
\end{proof}


\begin{lemma}\label{lem:lattice-near-point}
  For each point in $z\in\C$ there is a lattice point $u\in\Ztau$ not
  divisible by $\tau$ with
  \begin{equation*}
    \abs{z-u} < \tfrac12\!\abs{\tau} + 1 = \tfrac12 \sqrt{q} + 1
  \end{equation*}
\end{lemma}


\begin{proof}
  First note that if $\eta\in\Ztau$ is divisible by $\tau$, then $\eta+1$ and
  $\eta-1$ are not divisible by $\tau$. Consider the lines $z(z+\tau/2)$ and
  $z(z-\tau/2)$. One of these lines cuts a horizontal line with lattice points
  $\eta_k=j\tau+k$, for some fixed $j\in\Z$ and all $k\in\Z$, on it. This
  means that a lattice point $u$ can be found by first going from $z$ at
  most a distance of $\abs{\tau}/2$ and then at most $1$ on the horizontal
  line. Strictly smaller holds since both directions are linearly independent.
\end{proof}


\begin{proposition}\label{pro:compensate}
  If either $w\geq4$ and $q\geq11$ or $w\geq8$, then there is a possible
  compensate change.

  More precisely, if $w\geq5$ or $q\geq5$, then there is a digit $a$ such that
  $a+\tau^{w-1}$ is a digit as well. If either $w\geq4$ and $q\geq11$ or
  $w\geq8$, then there is a digit $b$ such that $b+\conj{\tau}\tau^{w-1}$ is a
  digit as well.
\end{proposition}


\begin{proof}
  The digit set $\cD$ is contained in $\tau^wV$. Consider the line form
  $\tau^{w-1}v_{-1/2}=\frac12\tau^{w-1}$ to
  $-\tau^{w-1}v_{-1/2}=-\frac12\tau^{w-1}$. From each end point of that line
  there is a lattice point not divisible by $\tau$ within a radius
  $\sqrt{q}/2+1$ by Lemma~\ref{lem:lattice-near-point}. With the quantities of
  Lemma~\ref{lem:distances-voronoi-lsd} one can easily check that the
  inequality
  \begin{equation*}
    \frac12 \sqrt{q}+1 < \abs{\tau}^w \frac{d_A}{\abs{\tau}} 
    = \sqrt{q}^{w-1} \left(\frac{\sqrt{q}}{2} - \frac{1}{4\sqrt{q}}\right)
  \end{equation*}
  holds for either $w\geq5$ or $q\geq5$ (fixing $w=2$ makes it easy to check
  the inequality for $q\geq5$, then use monotonicity in $w$; use the same
  argumentation starting with $q=2$ and $w\geq5$). Thus, the lattice points
  found above are in the interior of $\tau^wV$, and so are our desired digits
  $a$ and $a+\tau^{w-1}$. 

  Now we do similarly to get digits $b$ and $b+\conj{\tau}\tau^{w-1}$. We
  consider the line from $\tau^{w-1}v_{3/2}=-\frac{\conj{\tau}}{2}\tau^{w-1}$
  to $-\tau^{w-1}v_{3/2}=-\frac{\conj{\tau}}{2}\tau^{w-1}$ and have to check
  the inequality
  \begin{equation*}
    \frac12 \sqrt{q}+1 < \abs{\tau}^w \frac{d_B}{\abs{\tau}}
    = \sqrt{q}^{w-1} \frac{1}{4\sqrt{q}}.
  \end{equation*}
  That inequality is satisfied either if $w\geq4$ and $q\geq11$ or if $w\geq8$,
  which can again be checked easily by monotonicity arguments. For $w=2$ or
  $w=3$ the inequality is never satisfied.
\end{proof}


\section{Finding a non-optimal \wNAF{}}
\label{sec:strategy}


In contrast to Section~\ref{sec:non-opt-w}, where we have given an explicit
construction for a counterexample to minimality of the \wNAF{}, we use a
different approach here. It still builds up on the ideas mentioned in the
introduction of Part~\ref{sec:digits}, i.e., for our construction we consider
an element
\begin{equation*}
  z = A \tau^{w-1} + B \in \Ztau
\end{equation*}
with nonzero digits $A\in\cD$ and $B\in\cD$ with the following additional
properties. We want to find a change $\Delta$ with
\begin{itemize}
\item $\tau \divides (A+\Delta)$,
\item $\tau^2 \ndivides (A+\Delta)$,
\item $\tau^{-1}(A+\Delta) \not\in\cD$ and
\item $B - \Delta\tau^{w-1} \in \cD$.
\end{itemize}
We will restrict ourselves here to
\begin{equation*}
  \Delta\in\set{-1,1,\tau-p,-\tau+p},
\end{equation*}
which turns out to be a good choice. (Note that this restriction was already
used in Section~\ref{sec:small-change-Delta}; we only will relax it in
Proposition~\ref{pro:nonopt-q4w6})

Then the \wNAF{}-expansion of $C=\tau^{-1}(A+\Delta)$ has weight at
least~$2$, because it (its value) is not a digit. We obtain
\begin{equation*}
  z = A \tau^{w-1} + B 
  = C\tau^w + \left(B - \Delta\tau^{w-1}\right),
\end{equation*}
which shows that $z$ has a (non-\wNAF{}) expansion with weight~$2$ and a
\wNAF{}-expansion with weight at least~$3$. For all of our cases, the right
hand side of the previous equation can be rewritten in an expansion
\begin{equation*}
  z = E\tau^{2w} + F\tau^w + G
\end{equation*}
with some digits $E$, $F$ and $G=D=B - \Delta\tau^{w-1}$.

The finding of a point $A+\Delta$ is based on the main result of
Part~\ref{sec:dioph}, $A$ and $\Delta$ are discussed in
Section~\ref{sec:small-change-Delta}, and Section~\ref{sec:lsd} is
devoted to get digits $B$ and $G=D$. We just have to glue all the
results together, which is done in the proposition below and in the
next section. Alternatively, a direct search can be used to find those
lattice point configurations; we use this when Part~\ref{sec:dioph}
does not provide a result.


\begin{proposition}\label{pro:construct-ce-by-Rw}
  Suppose either $w\geq4$ and $q\geq11$ or $w\geq8$ (as in
  Proposition~\ref{pro:compensate}). Moreover, set
  \begin{equation*}
    C = \tau^{-1} \left( a\tau + b(q-p\tau) \right)
  \end{equation*}
  and suppose we have
  \begin{equation*}
    \tau C \in R_w
  \end{equation*}
  for some $a\in\Z$ with $q\ndivides a$ and $b\in\Z$. Then there exist digits $A$, $B$
  and $D$ such that
  \begin{equation*}
    A \tau^{w-1} + B = C \tau^w + D.
  \end{equation*}
\end{proposition}


Note that $C$ cannot be a digit because of the following reasons. The point
$\tau C$ lies in the rectangle~$R_w$, and thus $C$ is outside of $\tau^w V$, see
also Remark~\ref{rem:Rw-divide}. Since $C$ not divisible by $\tau$ (because
$q\ndivides \ell$) either, it has an expansion of weight at least~$2$. Using
this expansion leads to our desired counterexample.


\begin{proof}
  Consider the lattice point $\tau C \in R_w$. By
  Lemma~\ref{lem:point-in-int-V} there exists a
  $\Delta\in\set{-1,1,\tau-p,-\tau+p}$ such that
  \begin{equation*}
    A = \tau C - \Delta
  \end{equation*}
  lies in the interior of the scaled Voronoi cell $\tau^w V$. Since $C$ is a
  lattice point (and $\Delta$ is not divisible by $\tau$), the lattice point
  $A$ is not divisible by $\tau$. Therefore $A$ is a digit
  (see Section~\ref{sec:digit-sets}).

  Proposition~\ref{pro:compensate} gives us the digits $B$ and $D = B - \Delta
  \tau^{w-1}$. This completes the proof.
\end{proof}


\section{Collecting all Results}
\label{sec:proof-main}


In this final section, we prove the Theorems~\ref{thm:non-opt-large-w}
and~\ref{thm:main-non-opt} and Proposition~\ref{pro:main-algo}.


\begin{proof}[Proof of Theorem~\ref{thm:non-opt-large-w}]
  If we have a lattice point $\tau C \in R_w$ not divisible by $\tau^2$, then
  we are able to construct a counterexample by
  Proposition~\ref{pro:construct-ce-by-Rw}.
  Fortunately Proposition~\ref{pro:expl-comp-bound-Rw} provides
  the explicitly computable bound~$w_q$ and that for all $w\geq w_q$ a $\tau C$
  as above exists.
\end{proof}


Before proving Theorem~\ref{thm:main-non-opt} we need to consider one
special case first.


\begin{proposition}\label{pro:nonopt-q4w6}
  Let $q=4$ and $p\in\set{-1,1}$ and $w=6$. Set
  \begin{align*}
    & &
    E &= 1, \\
    A &= 7\tau - 66p, &
    F &= -16p\tau + 10, \\
    B &= -65, &
    G &= 10p\tau - 9.
  \end{align*}
  Then
  \begin{equation}\label{eq:nonopt-q4w6}
    A\tau^{5} + B = E\tau^{12} + F\tau^6 + G
  \end{equation}
  and both sides of the equation are valid digit expansions, i.e., the $6$-NAF
  is not a minimal digit expansion.
\end{proposition}

\begin{proof}
  We use a direct search following the same ideas as presented above
  (especially in Section~\ref{sec:strategy}). However, we have to
  relax our conditions on $\Delta$, in particular, we use
  $\Delta=-2p$. Moreover, we have $C=17p\tau-10$ as intermediate
  result in the construction.
\end{proof}


\begin{proof}[Proof of Proposition~\ref{pro:main-algo}
  and Theorem~\ref{thm:main-non-opt}]
  We start as in the proof of Theorem~\ref{thm:non-opt-large-w}, i.e., we need
  a lattice point $\tau C \in R_w$ not divisible by $\tau^2$.
  The main result of Part~\ref{sec:dioph}, namely
  Proposition~\ref{pro:find-lattice-point}, provides the existence of such a
  lattice point for a fixed~$q$ with finitely many (only a few)
  exceptional values for $w$. For these exceptions, we perform a
  direct search over all possible lattice points to get $\tau C$ (not
  lying inside~$R_w$) and construct the counterexample as described in
  Section~\ref{sec:strategy}. Note that a possible compensate change
  (Section~\ref{sec:lsd}) can be found by a lattice search as well.

  This construction of the actual counterexample and lattice search
  for the exceptions extends Algorithm~\ref{alg:test-fixed-q}; thus
  completes the proof of Proposition~\ref{pro:main-algo}.

  Applying this algorithm, the existence of counterexamples for all $w\geq4$ is
  shown; the only exceptions are $q=4$, $p\in\set{-1,1}$ and $w=6$
  which are handled separately by Proposition~\ref{pro:nonopt-q4w6}.

  Non-minimality of the cases $w=2$ and $w=3$ (and
  arbitrary $q\geq3$) is proven in
  Section~\ref{sec:non-opt-w}. Minimality for $q=2$ and
  $w\in\set{2,3}$ is shown
  in~\cite{Avanzi-Heuberger-Prodinger:2006:minim-hammin,
    Avanzi-Heuberger-Prodinger:2006:scalar-multip-koblit-curves,
    Gordon:1998}. This finishes the proof of Theorem~\ref{thm:main-non-opt}.
\end{proof}


\renewcommand{\MR}[1]{}

\bibliographystyle{amsplainurl}
\bibliography{koblitz2-non-optimal}

\providecommand{\Submitted}{Submitted} \providecommand{\availableat}{ available
  at } \providecommand{\alsoavailableat}{ also available at }
  \providecommand{\evavailableat}{earlier version available at }
  \providecommand{\toappearin}{To appear in } \providecommand{\toappear}{to
  appear} \providecommand{\inpreparation}{in preparation}
  \providecommand{\doi}[1]{\href{http://dx.doi.org/#1}{\path{doi:#1}}}
  \providecommand{\etc}{\emph{etc.}}\def\cprime{$'$}
\providecommand{\bysame}{\leavevmode\hbox to3em{\hrulefill}\thinspace}
\providecommand{\MR}{\relax\ifhmode\unskip\space\fi MR }
\providecommand{\MRhref}[2]{%
  \href{http://www.ams.org/mathscinet-getitem?mr=#1}{#2}
}
\providecommand{\href}[2]{#2}
\begin{thebibliography}{10}

\bibitem{avanzi:mywnaf}
Roberto Avanzi, \href{http://dx.doi.org/10.1007/978-3-540-30564-4_9}{\emph{A
  note on the signed sliding window integer recoding and a left-to-right
  analogue}}, {Selected Areas in Cryptography: 11th International Workshop, SAC
  2004, Waterloo, Canada, August 9-10, 2004, Revised Selected Papers}
  (H.~Handschuh and A.~Hasan, eds.), Lecture Notes in Comput. Sci., vol. 3357,
  Springer-Verlag, Berlin, 2005, pp.~130--143.

\bibitem{Avanzi-Heuberger-Prodinger:2006:minim-hammin}
Roberto Avanzi, Clemens Heuberger, and Helmut Prodinger,
  \href{http://dx.doi.org/10.1007/11693383_23}{\emph{Minimality of the
  {H}amming weight of the $\tau$-{NAF} for {K}oblitz curves and improved
  combination with point halving}}, Selected Areas in Cryptography: 12th
  International Workshop, SAC 2005, Kingston, ON, Canada, August 11--12, 2005,
  Revised Selected Papers (B.~Preneel and S.~Tavares, eds.), Lecture Notes in
  Comput. Sci., vol. 3897, Springer, Berlin, 2006, pp.~332--344. \MR{2241647
  (2007f:94028)}

\bibitem{Avanzi-Heuberger-Prodinger:2006:scalar-multip-koblit-curves}
\bysame, \href{http://dx.doi.org/10.1007/s00453-006-0105-9}{\emph{Scalar
  multiplication on {K}oblitz curves. {U}sing the {F}robenius endomorphism and
  its combination with point halving: {E}xtensions and mathematical analysis}},
  Algorithmica \textbf{46} (2006), 249--270. \MR{2291956 (2008a:94180)}

\bibitem{Avanzi-Heuberger-Prodinger:2011:redun-expan-i}
\bysame, \href{http://dx.doi.org/10.1007/s10623-010-9396-6}{\emph{Redundant
  $\tau$-adic expansions {I}: Non-adjacent digit sets and their applications to
  scalar multiplication}}, Des. Codes Cryptogr. \textbf{58} (2011), 173--202.
  \MR{2770310 (2012f:11010)}

\bibitem{Baker-Davenport:1969}
Alan Baker and Harold Davenport, \emph{The equations {$3x^2-2=y^2$} and
  {$8x^2-7=z^2$}}, Quart. J.~Math. Oxford \textbf{20} (1969), 129--137.
  \MR{MR0248079 (40 \#1333)}

\bibitem{Blake-Seroussi-Smart:1999}
Ian~F. Blake, Gadiel Seroussi, and Nigel~P. Smart, \emph{Elliptic curves in
  cryptography}, London Mathematical Society Lecture Note Series, vol. 265,
  Cambridge University Press, 1999.

\bibitem{Cassels:1959:GN}
John W.~S. Cassels, \emph{An introduction to the geometry of numbers}, Die
  Grundlehren der mathematischen Wissenschaften in Einzeldarstellungen mit
  besonderer Ber\"ucksichtigung der Anwendungsgebiete, Bd. 99 Springer-Verlag,
  Berlin-G\"ottingen-Heidelberg, 1959. \MR{0157947 (28 \#1175)}

\bibitem{Gordon:1998}
Daniel~M. Gordon, \emph{A survey of fast exponentiation methods}, J. Algorithms
  \textbf{27} (1998), 129--146. \MR{99g:94014}

\bibitem{Heuberger:2010:nonoptimality}
Clemens Heuberger,
  \href{http://dx.doi.org/10.1007/s11786-009-0014-9}{\emph{Redundant
  $\tau$-adic expansions {II}: {N}on-optimality and chaotic behaviour}}, Math.
  Comput. Sci. \textbf{3} (2010), 141--157. \MR{2608292 (2011b:11012)}

\bibitem{Heuberger-Krenn:2013:wnaf-analysis}
Clemens Heuberger and Daniel Krenn,
  \href{http://dx.doi.org/10.1016/j.jnt.2012.08.029}{\emph{Analysis of
  width-$w$ non-adjacent forms to imaginary quadratic bases}}, J.~Number Theory
  \textbf{133} (2013), no.~5, 1752--1808. \MR{3007130}

\bibitem{Heuberger-Krenn:2013:general-wnaf}
\bysame, \href{http://dx.doi.org/10.1007/s10474-013-0303-2}{\emph{Existence and
  optimality of $w$-non-adjacent forms with an algebraic integer base}}, Acta
  Math. Hungar. \textbf{140} (2013), no.~1--2, 90--104. \MR{3123865}

\bibitem{Heuberger-Krenn:2013:wnafs-optimality}
\bysame, \href{http://dx.doi.org/10.5802/jtnb.840}{\emph{Optimality of the
  width-$w$ non-adjacent form: General characterisation and the case of
  imaginary quadratic bases}}, J.~Th\'eor. Nombres Bordeaux \textbf{25} (2013),
  no.~2, 353--386. \MR{3228312}

\bibitem{Jedwab-Mitchell:1989}
Jonathan Jedwab and Chris~J. Mitchell, \emph{{Minimum weight modified
  signed-digit representations and fast exponentiation}}, Electron. Lett.
  \textbf{25} (1989), 1171--1172.

\bibitem{Koblitz:1992:cm}
Neal Koblitz,
  \href{http://dx.doi.org/10.1007/3-540-46766-1_22}{\emph{C{M}-curves with good
  cryptographic properties}}, Advances in cryptology---CRYPTO '91 (Santa
  Barbara, CA, 1991) (J.~Feigenbaum, ed.), Lecture Notes in Comput. Sci., vol.
  576, Springer, Berlin, 1992, pp.~279--287. \MR{94e:11134}

\bibitem{Koblitz:1998:ellip-curve}
\bysame, \emph{An elliptic curve implementation of the finite field digital
  signature algorithm}, Advances in cryptology---CRYPTO '98 (Santa Barbara, CA,
  1998), Lecture Notes in Comput. Sci., vol. 1462, Springer, Berlin, 1998,
  pp.~327--337. \MR{MR1670960 (99j:94052)}

\bibitem{Matveev:2000:explicit-lower-bound}
Eugene~M. Matveev, \emph{An explicit lower bound for a homogeneous rational
  linear form in logarithms of algebraic numbers. {II}}, Izv. Ross. Akad. Nauk
  Ser. Mat. \textbf{64} (2000), no.~6, 125--180. \MR{2002e:11091}

\bibitem{Meier-Staffelbach:1993:effic}
Willi Meier and Othmar Staffelbach, \emph{Efficient multiplication on certain
  nonsupersingular elliptic curves}, Advances in cryptology---CRYPTO '92 (Santa
  Barbara, CA, 1992) (Ernest~F. Brickell, ed.), Lecture Notes in Comput. Sci.,
  vol. 740, Springer, Berlin, 1993, pp.~333--344. \MR{MR1287863 (95e:94037)}

\bibitem{Miyaji-Ono-Cohen:1997:effic}
Atsuko Miyaji, Takatoshi Ono, and Henri Cohen, \emph{Efficient elliptic curve
  exponentiation}, Information and communications security. 1st international
  conference, ICICS '97, Beijing, China, November 11--14, 1997. Proceedings
  (Yongfei Han, Tatsuaki Okamoto, and Sihan Qing, eds.), Lecture Notes in
  Comput. Sci., vol. 1334, Springer-Verlag, 1997, pp.~282--290.

\bibitem{muirstinson:minimality}
James~A. Muir and Douglas~R. Stinson,
  \href{http://dx.doi.org/10.1090/S0025-5718-05-01769-2}{\emph{Minimality and
  other properties of the width-{$w$} nonadjacent form}}, Math. Comp.
  \textbf{75} (2006), 369--384.

\bibitem{Pethoe:1986}
Attila Peth{\"o} and Benjamin M.~M. de~Weger,
  \href{http://dx.doi.org/10.2307/2008185}{\emph{Products of prime powers in
  binary recurrence sequences. {I}. {T}he hyperbolic case, with an application
  to the generalized {R}amanujan-{N}agell equation}}, Math. Comp. \textbf{47}
  (1986), no.~176, 713--727. \MR{856715 (87m:11027a)}

\bibitem{Phillips-Burgess:2004:minim-weigh}
Braden Phillips and Neil Burgess,
  \href{http://dx.doi.org/10.1109/TC.2004.14}{\emph{Minimal weight digit set
  conversions}}, IEEE Trans. Comput. \textbf{53} (2004), 666--677.

\bibitem{Reitwiesner:1960}
George~W. Reitwiesner,
  \href{http://dx.doi.org/10.1016/S0065-2458(08)60610-5}{\emph{Binary
  arithmetic}}, Advances in Computers, {V}ol. 1, Academic Press, New York,
  1960, pp.~231--308. \MR{0122018 (22 \#12745)}

\bibitem{Ribenboim:2001}
Paulo Ribenboim, \emph{Classical theory of algebraic numbers}, Universitext,
  Springer-Verlag, New York, 2001. \MR{1821363 (2002e:11001)}

\bibitem{SageMath:2016:7.0}
The~SageMath Developers, \emph{{SageMath} {M}athematics {S}oftware ({V}ersion
  7.0)}, 2016, \url{http://www.sagemath.org}.

\bibitem{Solinas:1997:improved-algorithm}
Jerome~A. Solinas, \emph{An improved algorithm for arithmetic on a family of
  elliptic curves}, Advances in Cryptology --- {CRYPTO} '97. 17th annual
  international cryptology conference. Santa Barbara, {CA}, {USA}. August
  17--21, 1997. Proceedings (B.~S. Kaliski, jun., ed.), Lecture Notes in
  Comput. Sci., vol. 1294, Springer, Berlin, 1997, pp.~357--371.

\bibitem{Solinas:2000:effic-koblit}
\bysame, \href{http://dx.doi.org/10.1023/A:1008306223194}{\emph{Efficient
  arithmetic on {K}oblitz curves}}, Des. Codes Cryptogr. \textbf{19} (2000),
  195--249. \MR{2002k:14039}

\end{thebibliography}


\end{document}